\RequirePackage{plautopatch}
\documentclass[smallextended,natbib]{svjour3}
\smartqed  
\usepackage{amsmath,amssymb}
\usepackage{graphicx}
\usepackage{url}
\usepackage{mathtools}\mathtoolsset{showonlyrefs=false}
\usepackage{algpseudocode,algorithm}
\usepackage{subfig}

\usepackage{helvet,courier}
\usepackage{newtxtext,newtxmath}

\newcommand{\refalg}[1]{Algorithm~\ref{#1}}

\newcommand{\reffig}[1]{Figure~\ref{#1}}

\newcommand{\finbox}{\nolinebreak\hfill{\small $\blacksquare$}}

\newcommand{\MIN}{\mathop{\mathrm{Minimize}}}

\newcommand{\ST}{\mathop{\mathrm{subject~to}}}

\renewcommand{\Re}{\ensuremath{\mathbb{R}}}

\newcommand{\bi}[1]{\ensuremath{\boldsymbol{#1}}}

\newcommand{\rr}[1]{\ensuremath{\mathrm{#1}}}

\begin{document}

\title{%
  Data-driven confidence bound for structural response 
  using segmented least squares: 
  a mixed-integer programming approach
  }

\author{%
  Yoshihiro Kanno
  }

\institute{
  Yoshihiro Kanno \at
    Corresponding author. 
    Address: 
    Mathematics and Informatics Center, 
    The University of Tokyo, 
    Hongo 7-3-1, Tokyo 113-8656, Japan.
    E-mail: \texttt{kanno@mist.i.u-tokyo.ac.jp}. 
  }

\date{Received: date / Accepted: date}

\maketitle

\begin{abstract}
  As one of data-driven approaches to computational mechanics in 
  elasticity, this paper presents a method finding a bound for 
  structural response, taking uncertainty in a material data set into 
  account. 
  For construction of an uncertainty set, we adopt the segmented least 
  squares so that a data set that is  not fitted well by the linear 
  regression model can be dealt with. 
  Since the obtained uncertainty set is nonconvex, the optimization 
  problem solved for the uncertainty analysis is nonconvex. 
  We recast this optimization problem as a mixed-integer programming 
  problem to find a global optimal solution. 
  This global optimality, together with a fundamental property of the 
  order statistics, guarantees that the obtained bound for the 
  structural response is conservative, in the sense that, at least a 
  specified confidence level, probability that the structural response 
  is in this bound is no smaller than a specified target value. 
  We present numerical examples for three different types of 
  skeletal structures. 
\end{abstract}



\section{Introduction}

Conventional computational mechanics assumes a model of the constitutive 
law specific to a material. 
For example, in the static equilibrium analysis of elastic solids and 
structures, the constitutive law relates the stress to the strain. 
It is usual that the constitutive law is determined through empirical 
modeling of the stress--strain relation, followed by calibration of the 
parameters in the model. 
In contrast, data-driven approaches to computational mechanics have 
recently attracted substantial attention, where a material data set is 
utilized directly for the structural analysis, without resorting to 
the conventional process of modeling and calibration. 

The contribution that activates emerging study of data-driven methods in 
computational mechanics is \citet{KO16}, where, provided that a data set 
of material experiments (i.e., a set of pairs of observed stress and 
strain values) is given, the structural response at the static 
equilibrium state is estimated. 
Specifically, based on the conventional finite element method, 
the stress and strain at each integration point are considered unknowns 
subjected to the compatibility condition and the force-balance equation. 
Then the pair of stress and strain values that has the minimum distance, 
in some sense, to the data set is declared as the equilibrium state. 
\citet{KO16} defined the distance from a point in the stress--strain 
space to a data set as the Euclidean distance (with a scaling) from the 
point to the closest data point. 
This methodology is sometimes called 
the {\em distance-minimizing method\/} \citep{NRK20,PBG23}, 
and has rapidly been applied to diverse problem settings, including 
elasticity with geometrical nonlinearity \citep{NK18}, 
history-dependent elastoplastic 
problems \citep{CH22,DLS24,EKRSO19,PBG23,PBS23}, 
brittle fracture \citep{CDSO20}, and 
diffusion problems \citep{NRK20}. 
It has also been extended to multi-scale modelings of, e.g., 
composites \citep{XYYHGBZBH20}, 
bone tissue \citep{MARDDDS20}, and granular materials \citep{KSOA21}. 
Moreover, instead of a data set of material experiments, 
usage of a data set of structural experiments has also been 
studied \citep{DCLV19,LCRSV18,SLO19}. 
It is worth noting that, for the multi-scale modeling~\citep{KSOA21}, 
etc., a material data set is sometimes generated by using 
numerical simulation. 
This is analogous to preparation of a data set by numerical material 
tests using the computational homogenization method 
to identify the macroscopic material 
properties of composites and polycrystalline metals with 
periodic micro-structures \citep{CSY12,TW07,TKHIY13,WT10}. 

An alternative methodology of data-driven computational mechanics 
adopts the notion of so-called {\em constitutive manifold\/}. 
Namely, the points satisfying the constitutive law usually do not exist 
ubiquitously in whole of the stress--strain space, but lie on a manifold 
with a lower dimension than the space. 
\citet{IAcAGCC18,IBAAcCLC17} proposed to use a locally linear embedding, 
which is one of popular methods for the manifold learning, to capture the 
constitutive manifold. 
\citet{HC20,HLLC21}, and \citet{SJC23} proposed methods that construct a 
local convex envelop from the neighboring data points. 
The author proposed kernel-based methods for extracting the constitutive 
manifold~\citep{Kan21,Kan21TAML}. 

Other approaches to data-driven computational mechanics are 
also found in literature. 
\citet{LP23} proposed a variant of the support vector regression that is 
robust against presence of outliers, and applied it to the static 
equilibrium analysis of an elastic truss. 
\citet{TLQYSMLG20,TYQFLG21,TZYLLG19} developed a method recovering the 
three-dimensional stress--strain relation from material data of the 
uniaxial experiments. 
This method was applied to continuum-based topology optimization in 
conjunction with the moving morphable void method~\citep{GDWMZSYTG23}. 
\citet{PRO23} supposed that the material response is random. 
A certain likelihood is given for each data point of a material data set, 
and the distribution of material response is estimated by minimizing the 
Kullback--Leibler divergence. 
Accordingly, \citet{PRO23} evaluated the expected value of the 
structural response. 

Motivated by the observation that, in the presence of even a single 
outlier, the method proposed by \citet{KO16} can possibly converge to 
a very improper solution (see \citep[Remark~2]{Kan19MIP} and 
\citep[Remark~1]{GDLT21}), 
\citet{GDLT21} proposed a distinctive method that introduces the concept 
of worst-case analysis under non-probabilistic 
uncertainty \citep{Kan20three} to the data-driven computational mechanics. 
Specifically, when a set of pairs of stress and strain values is given, 
this method considers an ellipsoid including all the data points as the 
uncertainty set. 
A bound of the structural response is then obtained by solving an 
optimization problem under the constraints that each pair of the element 
stresses and strains belongs to the ellipsoid, together with the 
compatibility condition and the force-balance equation. 
It is worth noting that this optimization problem is convex, and hence 
can be solved globally, which guarantees the conservativeness of the 
obtained response bound. 
Thus, the distinguished feature of the method of \citet{GDLT21} is that 
it provides upper and lower bounds, rather than a single value, of the 
quantity of interest. 
This is attractive from the perspective of 
the {\em uncertainty quantification\/} (UQ)~\citep{GHO17}, 
because the aleatory uncertainty, a.k.a.\ the natural variability, 
inevitably and intrinsically exists in material properties. 

Inspired by \citet{GDLT21}, \citet{Kan23} proposed a method finding a 
response bound with a specified confidence. 
More precisely, it is guaranteed that, at least a specified confidence 
level, the probability that the structural response belongs to the 
obtained bound is no smaller than the target reliability. 
Here, the material property is not considered deterministic but 
stochastic, and given material data are supposed to be 
independent and identically distributed samples 
drawn from a distribution. 
A key is that the number of data points which the uncertainty set should 
include to guarantee the specified confidence level can be given by a 
fundamental property of the order statistics. 
Accordingly, the method does not require any modeling of the 
distribution, and hence is viewed as a model-free and data-driven 
UQ method. 
It is worth noting that the notion of confidence level is borrowed from 
studies on the reliability-based design optimization under 
uncertainty in the input distribution 
\citep{HYYZWW22,IKK18,JCDL20,JCL19,Kan19,Kan20,MCCGLG18,WHYWG20}. 

A drawback of \citet{Kan23} is that, as it uses the linear regression 
to construct the uncertainty set, its applicability is limited to 
approximately linear material data: For a nonlinear elastic material, 
the bound obtained by this method can drastically overestimate the 
structural response. 
Use of the linear regression is motivated by the desire to guarantee the 
global optimality of solutions of the optimization problems for finding 
the response bound; a local optimal solution in general underestimates 
the structural response. 
In fact, the method of \citet{GDLT21} has the same drawback, 
where it solves an optimization of linear objective function under convex 
quadratic constraints, while the method of \citet{Kan23} solves a linear 
programming problem. 
To deal with this drawback, \citet{HLDTG23} proposed to use a local 
uncertainty set for each element stress and strain, rather than a unique 
global ellipsoid that is common to all the elements. 
Specifically, this method adopts 
the convex hull of data points only in the neighborhood of 
each incumbent element stress and strain point. 
As a result, the method lacks guarantee of the global optimality, and 
hence the obtained solution in general underestimates the structural 
response. 
In contrast, this study 
attempts to deal with the drawback mentioned above with 
maintaining guarantee of the global optimality. 
Specifically, we extend the method of \citet{Kan23} by adopting the 
segmented least squares, instead of the linear regression. 
Although the induced uncertainty set is nonconvex, we show that this set 
is handled within the framework of 
{\em mixed-integer linear programming\/} (MILP), 
which is a key to guaranteeing the global optimality. 

One of potential advantages of the proposed method, as well as other 
methods based on uncertainty analysis \citep{GDLT21,HLDTG23,Kan23}, 
over the distance-minimizing methods 
\citep{CDSO20,CH22,DCLV19,DLS24,EKRSO19,KSOA21,KO16,LCRSV18,MARDDDS20,NK18,NRK20,PBG23,PBS23,SLO19,XYYHGBZBH20} 
is that it does not require a large number of data points. 
For example, in the numerical experiments, 
\citet[section~3]{KO16}, \citet[section~4.2]{NK18}, 
and \citet[section~5.3]{PBS23} used data sets consisting of about 
$10^{6}$ data points; 
\citet[section~2.8.2]{NRK20} used up to about $1.6\times 10^{5}$ data 
points; and \citet[section~5]{HC20} used up to about $5.1\times 10^{5}$ 
data points. 
In contrast, the methods proposed in \cite{Kan23} and this paper can be 
applied even when the number of data points is relatively small 
(as demonstrated in section~\ref{sec:ex}), where the upper bounds for 
the confidence level and target reliability become small. 
Precisely, by putting $\tilde{p}=r$ in \eqref{eq:def.tilde.p} in 
section~\ref{sec:framework}, we see that for the specified target 
reliability $1-\epsilon \in ]0,1[$ the upper bound for the realizable 
confidence level is $1-\delta=1-(1-\epsilon)^{r}$, where $r$ is the 
number of data points. 
For example, if the data set has $r=500$ data points and we choose 
$\epsilon=0.01$, then the minimum value of $\delta$ is about 
$6.57\times 10^{-3} \simeq (1-\epsilon)^{r}$.

The remaining paper is organized as follows: 
Section~\ref{sec:framework} revisits the framework proposed in 
\cite{Kan23} to introduce some useful definitions, and 
clarifies contributions of this study. 
Section~\ref{sec:regression} presents a {\em mixed-integer programming\/}
(MIP) formulation for the segmented least squares, which is used to 
construct an uncertainty set in the stress--strain space. 
Section~\ref{sec:uncertainty} shows that we can express the uncertainty 
set as some linear inequalities with some 0-1 variables, which enables 
us to formulate the optimization problem for finding a bound for the 
quantity of interest as a MIP problem. 
For simplicity, a concrete MIP formulation for finding a bound for the 
structural response is presented only for trusses in 
section~\ref{sec:truss}. 
In section~\ref{sec:ex}, we perform numerical experiments on three 
types of skeletal structures. 
Section~\ref{sec:conclusions} presents the conclusions.

\section{Framework of computation with confidence}
\label{sec:framework}

For completeness and ease of comprehension, we first overview the 
framework of the method proposed in the previous work \citep{Kan23}. 
We next identify the aims and contributions of this study. 

\subsection{Confidence bound for structural response}
\label{sec:framework.bound}

Suppose that we attempt to find a bound for the static 
response of an elastic structure. 
Assume that the quantity of interest, denoted by $q$, depends on a random 
vector $\bi{z} \in \Re^{c}$, where $c$ is the number of random variables. 
As a forecast of the structural response, consider a bound 
such that the probability that the realization of $q(\bi{z})$ is not 
included in this bound is no greater than the specified 
value, $\epsilon \in ]0,1[$, where $1-\epsilon$ is called 
the {\em target reliability\/}. 
Note that, in this paper, we use $]a,b[$ and $[a,b]$ to 
denote the open interval and closed interval, respectively, 
between $a\in\Re$ and $b\in\Re$. 
If we know the joint distribution function, $F$, of $\bi{z}$, then 
we seek to find $\underline{q} \in \Re$ and $\overline{q} \in \Re$, 
i.e., lower and upper bounds for the quantity of interest, respectively, 
satisfying 
\begin{align}
  \rr{P}\{ q(\bi{z}) \in [\underline{q},\overline{q}] \} \ge 1-\epsilon . 
  \label{eq:confidence.inequality.2}
\end{align}

In reality, we have no knowledge on $F$. 
Instead, suppose that we are given a data set of $\bi{z}$, 
which consists of a finite number of 
continuous independent and identically distributed samples, 
denoted by $\check{\bi{z}}_{1},\dots,\check{\bi{z}}_{r}\in\Re^{c}$, 
drawn from $F$.
We consider a lower bound constraint for the probability that 
constraint \eqref{eq:confidence.inequality.2} is satisfied, i.e., 
\begin{align}
  \rr{P}_{F}\bigl\{ 
  \rr{P}\{ q(\bi{z}) \in [\underline{q},\overline{q}] \} \ge 1-\epsilon 
  \bigr\} \ge 1-\delta , 
  \label{eq:confidence.inequality.1}
\end{align}
where $\delta \in ]0,1[$ is a specified value, 
$1-\delta$ is called the {\em target confidence level\/}, 
and $\rr{P}_{F}\{ \,\cdot\, \}$ means the probability 
taken for all possible $F$ for which 
$\check{\bi{z}}_{1},\dots,\check{\bi{z}}_{r}$ are 
continuous independent and identically distributed samples. 

The data-driven concept in \citep{Kan23} is based on the order 
statistics, and does not resort to any empirical modeling of $F$. 
Let $\tilde{p}$ denote the minimum integer satisfying 
\begin{align}
  \sum_{k=\tilde{p}}^{r} 
  \binom{r}{k} (1-\epsilon)^{k} \epsilon^{r-k} \le \delta . 
  \label{eq:def.tilde.p}
\end{align}
We use $Z \subset \Re^{c}$ to denote a set containing $\tilde{p}$ 
samples among given samples 
$\check{\bi{z}}_{1},\dots,\check{\bi{z}}_{r}$. 
It follows from Theorem~2.1 and the discussion in section~3 in 
\citet{Kan23} that $Z$ satisfies 
\begin{align*}
  \rr{P}_{F}\bigl\{ 
  \rr{P}\{ \bi{z} \in Z \} \ge 1-\epsilon 
  \bigr\} \ge 1-\delta . 
\end{align*}
Accordingly, we see that $\underline{q}$ and $\overline{q}$ defined by 
\begin{align}
  \underline{q} 
  &= \min \{ q(\bi{z}) \mid \bi{z} \in Z \} , 
  \label{eq:def.q_min} \\
  \overline{q} 
  &= \max \{ q(\bi{z}) \mid \bi{z} \in Z \} 
  \label{eq:def.q_max}
\end{align}
satisfy \eqref{eq:confidence.inequality.1}.

\subsection{Contributions of this study}

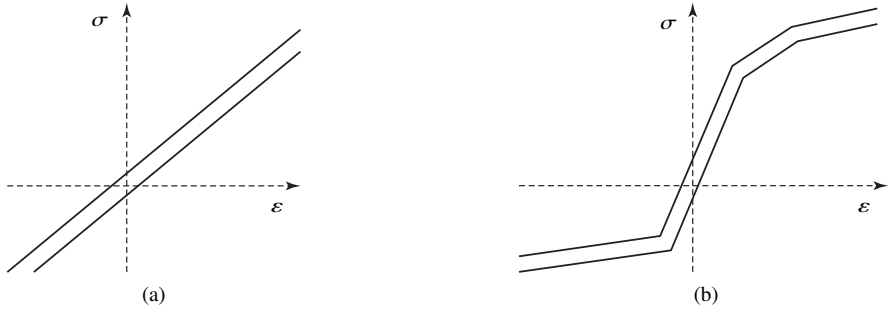
\begin{figure}[tbp]
  \centering
  \subfloat[]{
  \label{fig:set_0}
  \includegraphics[scale=0.45]{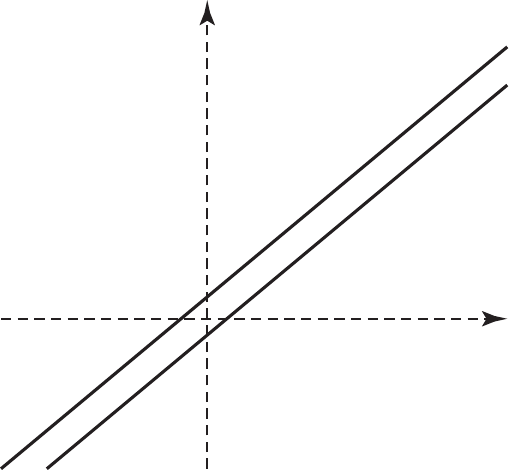}
  \begin{picture}(0,0)
    \put(-170,-50){
    \put(87,143){ $\sigma$ }
    \put(154,73){ $\varepsilon$ }
    }
  \end{picture}
  }
  \hfill
  \subfloat[]{
  \label{fig:set_1}
  \includegraphics[scale=0.45]{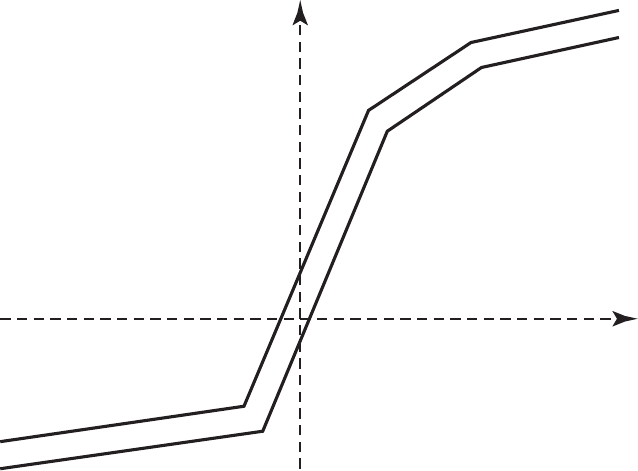}
  \begin{picture}(0,0)
    \put(-170,-50){
    \put(80,142){ $\sigma$ }
    \put(154,73){ $\varepsilon$ }
    }
  \end{picture}
  }
  \caption[]{Uncertainty set $C(\tau)$ of $(\varepsilon,\sigma)$. 
  \subref{fig:set_0}~The one considered in \cite{Kan23}; and 
  \subref{fig:set_1}~an example handled in this study. }
  \label{fig:set_0_1}
\end{figure}

As a concrete example consider a truss, where the constitutive law 
relates the uniaxial stress $\sigma\in \Re$ to the uniaxial 
strain $\varepsilon\in\Re$. 
We use $\sigma_{e}\in\Re$ and $\varepsilon_{e}\in\Re$ to 
denote the stress and strain of member $e$ $(e=1,\dots,m)$, respectively, 
where $m$ is the number of members of the truss. 

\begin{figure}[bp]
  \centering
  \subfloat[]{
  \label{fig:for_figure_prog3_tri_data_set}
  \includegraphics[scale=0.50]{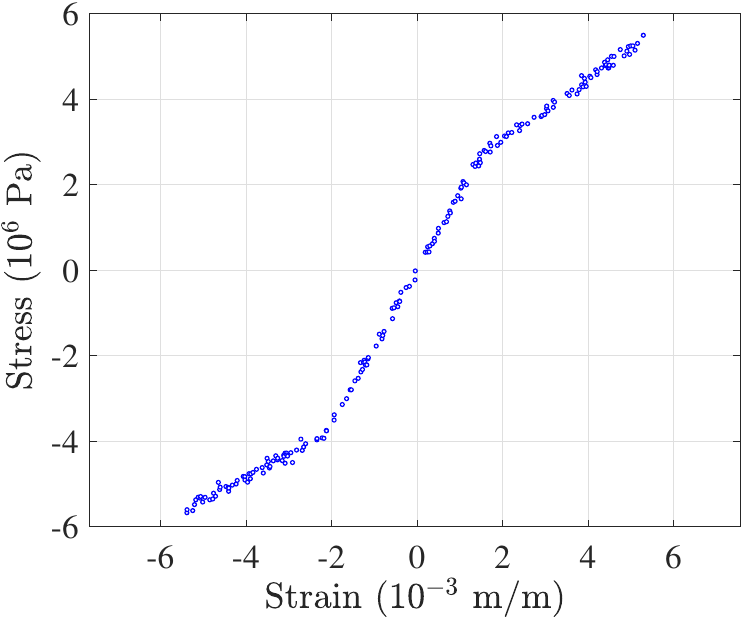}
  }
  \par
  \subfloat[]{
  \label{fig:for_figure_prog3_segmented}
  \includegraphics[scale=0.50]{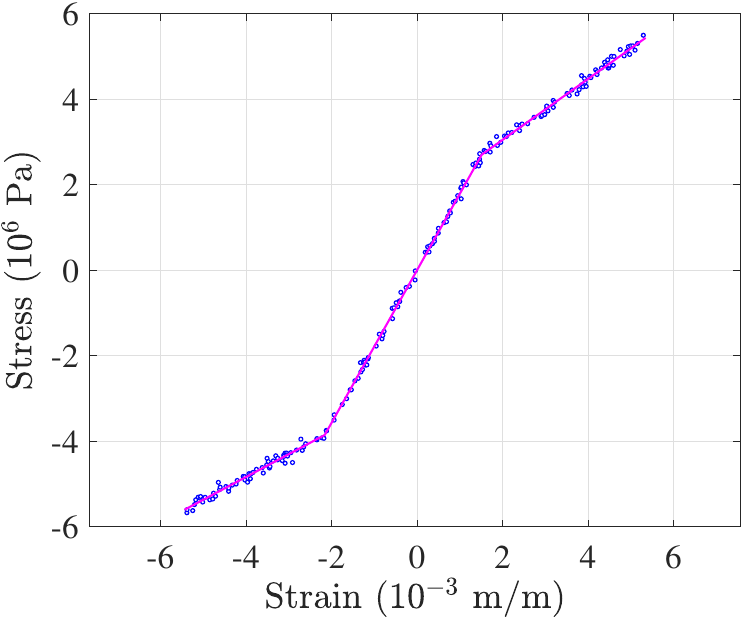}
  }
  \hfill
  \subfloat[]{
  \label{fig:for_figure_prog3_seg_bound}
  \includegraphics[scale=0.50]{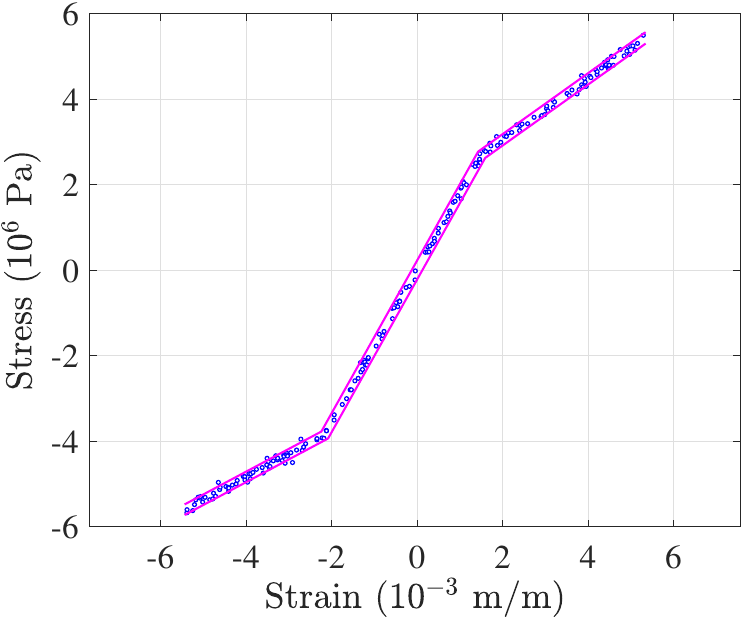}
  }
  \caption[]{Procedure of construction of $C(\tau)$ 
  proposed in this study. 
  \subref{fig:for_figure_prog3_tri_data_set}~To a given data set 
  \subref{fig:for_figure_prog3_segmented}~we first apply the segmented 
  least squares; 
  \subref{fig:for_figure_prog3_seg_bound}~and next we obtain $C(\tau)$ by 
  determining the value of $\tau$. }
\end{figure}

Suppose that the constitutive law essentially has uncertainty. 
We are given a set of continuous independent and identically distributed 
samples, denoted by 
$D = \{ (\check{\varepsilon}_{1},\check{\sigma}_{1}),\dots, 
  (\check{\varepsilon}_{r},\check{\sigma}_{r}) \}$, 
drawn from the probability distribution that $(\varepsilon,\sigma)$ 
follows. 
We use $C(\tau) \subset \Re^{2}$ to denote 
the uncertainty set of $(\varepsilon,\sigma)$, where $\tau>0$ is a 
parameter representing the magnitude of uncertainty. 
In \cite{Kan23}, we adopted a set depicted in \reffig{fig:set_0} 
as $C(\tau)$, i.e., 
\begin{align}
  C(\tau) 
  = \{ (\varepsilon,\sigma) \in \Re \times \Re
  \mid
  | \alpha \varepsilon + \beta \sigma - \gamma | \le \tau   \} , 
  \label{eq:linear.envelope}
\end{align}
where $\alpha$, $\beta$, $\gamma \in \Re$ are constants. 
We set the value of $\tau$ so that $C(\tau)$ includes $\tilde{p}$ 
samples among $(\check{\varepsilon}_{1},\check{\sigma}_{1}),\dots, 
  (\check{\varepsilon}_{r},\check{\sigma}_{r})$. 
Then, from the discussion in 
section~\ref{sec:framework.bound} we obtain 
\begin{align*}
  \rr{P}_{F}\bigl\{ 
  \rr{P}\{ (\varepsilon,\sigma) \in C(\tau) \} \ge 1-\epsilon 
  \bigr\} \ge 1-\delta . 
\end{align*}
Accordingly, we can formulate an optimization problem corresponding to 
problem \eqref{eq:def.q_min} as follows: 
The optimization variables are the nodal displacements, member stresses, 
and member strains. 
The constraints are (i) the compatibility relations between the nodal 
displacements and the member strains, 
(ii) the force-balance equations in terms of the member stresses and the 
nodal external forces, and (iii) the inclusions 
\begin{align*}
  (\varepsilon_{e}, \sigma_{e}) \in C(\tau) , 
  \quad e=1,\dots,m . 
\end{align*}
With this setting we minimize the quantity of interest to obtain 
its lower bound $\underline{q}$ satisfying \eqref{eq:confidence.inequality.1}. 
Also, an upper bound $\overline{q}$ satisfying 
\eqref{eq:confidence.inequality.1} is obtained by maximizing the  
quantity of interest \cite{Kan23}. 

The aim of this study is to develop a method that can adopt a set 
depicted in \reffig{fig:set_1} as $C(\tau)$. 
Suppose that the data set shown in 
\reffig{fig:for_figure_prog3_tri_data_set} is given as a set of samples 
of $(\varepsilon,\sigma)$. 
If we adopt the set shown in \reffig{fig:set_0} as $C(\tau)$, 
then the solution obtained by the method in \citep{Kan23} 
becomes too conservative, i.e., the solution extremely overestimates 
the structural response. 
In such a case the set shown in \reffig{fig:set_1} mitigates 
overestimate drastically, which is exactly the motivation of this study. 
The key is that this paper presents a method to obtain the global 
optimal solutions of the optimization problems in 
\eqref{eq:def.q_min} and \eqref{eq:def.q_max}. 
It is worth noting that a local optimal solution of 
\eqref{eq:def.q_min} or \eqref{eq:def.q_max} does not necessarily 
satisfies \eqref{eq:confidence.inequality.1}; the global optimality is 
crucial to ensure \eqref{eq:confidence.inequality.1}. 
If $C(\tau)$ is defined by \eqref{eq:linear.envelope} (i.e., 
\reffig{fig:set_0}), the optimization problems described above are 
formulated as {\em linear programming\/} problems. 
Therefore, in the previous work \citep{Kan23} it is straightforward to 
guarantee the global optimality of the proposed method. 
In contrast, the set in \reffig{fig:set_1} considered in this study is 
{\em nonconvex\/}. 
Therefore, global optimization is highly nontrivial. 
As a major contribution of this study, 
we show that the optimization problems described above can 
be formulated as {\em mixed-integer linear programming\/} (MILP) problems, 
which can be solved globally with, e.g., a branch-and-cut method. 

In this study, we suppose that a data set such as 
\reffig{fig:for_figure_prog3_tri_data_set} is given as 
the result of material experiments. 
We first find some lines that fit this data set as shown in 
\reffig{fig:for_figure_prog3_segmented}. 
We show that this segmented least squares can be formulated as 
a {\em mixed-integer second-order cone programming\/} (MISOCP) 
problem.\footnote{%
The problem is formulated as a mixed-integer programming problem with 
convex quadratic constraints, which can be recast as 
an MISOCP problem. }
Next, we construct uncertainty set $C(\tau)$ 
shown in \reffig{fig:for_figure_prog3_seg_bound}, and formulate 
an MILP problem to obtain a bound for the quantity of interest.

\section{Segmented least squares by mixed-integer programming}
\label{sec:regression}

This section presents segmented least squares problem that finds 
\reffig{fig:for_figure_prog3_segmented} when the data set in 
\reffig{fig:for_figure_prog3_tri_data_set} is given. 
Section~\ref{sec:regression.problem} elucidates the problem setting. 
Section~\ref{sec:regression.integer} formulates the problem as an MISOCP 
problem. 

\subsection{Problem setting}
\label{sec:regression.problem}

\begin{figure}[tbp]
  \centering
  \includegraphics[scale=1.00]{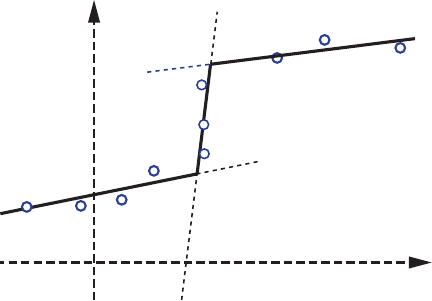}
  \begin{picture}(0,0)
    \put(-228,-55){
    \put(214,63){{\small $\varepsilon$}}
    \put(51,188){{\small $\sigma$}}
    \put(27,90){{\small $\check{\varepsilon}_{1}$}}
    \put(88,122){{\small $\check{\varepsilon}_{4}$}}
    \put(120,123){{\small $\check{\varepsilon}_{5}$}}
    \put(103,156){{\small $\check{\varepsilon}_{7}$}}
    \put(147,162){{\small $\check{\varepsilon}_{8}$}}
    \put(206,167){{\small $\check{\varepsilon}_{10}$}}
    \put(65,90){{\small $D_{1}$}}
    \put(119,144){{\small $D_{2}$}}
    \put(167,164){{\small $D_{3}$}}
    }
  \end{picture}
  \caption{Segmented least squares with three straight lines. 
  ``$\circ$'' denotes the data point 
  $(\check{\varepsilon}_{l},\check{\sigma}_{l})$ $(l=1,\dots,10)$ 
  numbered as $\check{\varepsilon}_{1} < \check{\varepsilon}_{2} 
  < \dots < \check{\varepsilon}_{10}$. }
  \label{fig:segment_4_rev}
\end{figure}
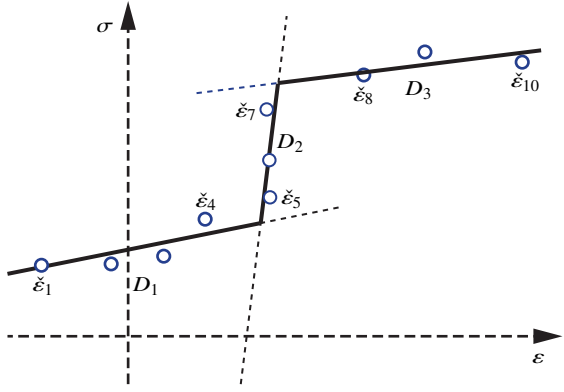

Let 
$D = \{ (\check{\varepsilon}_{1},\check{\sigma}_{1}),\dots, 
  (\check{\varepsilon}_{r},\check{\sigma}_{r})\}$ 
denote the data set consisting of pairs of the observed uniaxial strain 
and stress values. 
Without loss of generality, assume that the data points are numbered 
so that inequalities 
\begin{align*}
  \check{\varepsilon}_{1} < \check{\varepsilon}_{2} 
  < \dots < \check{\varepsilon}_{r}
\end{align*}
are satisfied. 
We attempt to find at most $k$ straight lines that best fit 
this data set; \reffig{fig:segment_4_rev} shows an example 
such that the data points approximately lie on either of the 
three straight lines. 
For each $i=1,\dots,k$, define line $\ell_{i}$ by 
\begin{align}
  \ell_{i} : 
  \quad
  \alpha_{i} \varepsilon + \beta_{i} \sigma = \gamma_{i} , 
\end{align}
where $\alpha_{i}$, $\beta_{i}$, $\gamma_{i} \in \Re$ are parameters 
that are to be determined.

This problem, the segmented least squares \citep[section~6.3]{KT06}, 
can be stated formally as follows. 
Let $\{ D_{1},\dots,D_{k} \}$ be a partition of $D$, i.e., 
$D_{1},\dots,D_{k} \subseteq D$, $D_{1} \cup \dots \cup D_{k} = D$, 
and $D_{i} \cap D_{i'} = \emptyset$ for $i \not= i'$. 
Some of $D_{1},\dots,D_{k}$ can possibly be empty, 
which means that the number of straight lines used for data fitting can 
be less than $k$. 
The data points belonging to nonempty $D_{i}$ have consecutive 
subscripts: 
In other words, with some integers 
$l_{1}$ and $l_{2}$ $(1 \le l_{1} \le l_{2} \le r)$, we can write 
\begin{align*}
  D_{i} = \{
  (\check{\varepsilon}_{l_{1}},\check{\sigma}_{l_{1}}),
  (\check{\varepsilon}_{l_{1}+1},\check{\sigma}_{l_{1}+1}),\dots, 
  (\check{\varepsilon}_{l_{2}-1},\check{\sigma}_{l_{2}-1}),
  (\check{\varepsilon}_{l_{2}},\check{\sigma}_{l_{2}}) \} . 
\end{align*}
Then we attempt fit line $\ell_{i}$ to the data points in $D_{i}$. 
If we adopt the least squares, then the error of the regression model 
for data point 
$(\check{\varepsilon}_{l},\check{\sigma}_{l}) \in D_{i}$ is given by 
$(\alpha_{i}\check{\varepsilon}_{l} 
  + \beta_{i}\check{\sigma}_{l} - \gamma_{i})^{2}$. 
As the number of nonempty $D_{i}$'s increases (i.e., as 
the number of breakpoints of the regression model increases), 
the sum of the errors obviously decreases. 
Therefore, we attempt to find $D_{1},\dots,D_{k}$ and 
$\alpha_{i}$, $\beta_{i}$, $\gamma_{i}$ $(i=1,\dots,k)$ that minimize 
the sum of (i) the number of nonempty $D_{i}$'s 
multiplied by penalty parameter $\mu > 0$ 
and (ii) the sum of squared errors of $\ell_{i}$ from 
$(\check{\varepsilon}_{l},\check{\sigma}_{l}) \in D_{i}$ $(i=1,\dots,k)$. 
If $\mu$ is sufficiently large, then at the optimal solution we have 
$D_{2}=\dots=D_{k}=\emptyset$ and hence our problem coincides with the 
conventional linear regression by the least squares. 
Alternatively, when $\mu>0$ is small, the number of breakpoints of the 
regression model increases.

\subsection{Mixed-integer programming formulation}
\label{sec:regression.integer}

This section presents a MIP formulation for the 
segmented least squares stated in section~\ref{sec:regression.problem}. 

In the following formulation, we introduce 0-1 variables 
representing which one of $\ell_{1},\dots,\ell_{k}$ is fitted to 
each data point $(\check{\varepsilon}_{l},\check{\sigma}_{l})$ 
$(l=1,\dots,r)$. 
Namely, for $(\check{\varepsilon}_{l},\check{\sigma}_{l})$, 
we set the values of $t_{li} \in \{ 0,1 \}$ $(i=1,\dots,k)$ as follows: 
\begin{align}
  (\check{\varepsilon}_{l},\check{\sigma}_{l}) \in D_{i}
  \quad\Leftrightarrow\quad
  \begin{dcases*}
    t_{l1} = \dots = t_{li} = 1 , \\
    t_{l,i+1} = \dots = t_{lk} = 0 .
  \end{dcases*}
  \label{eq:def.binary.t.1}
\end{align}
It is natural to fix 
\begin{align}
  t_{1,1}=1 , 
  \quad
  t_{1,2}=0 , 
  \label{eq:constraint.for.D1}
\end{align}
i.e., $(\check{\varepsilon}_{1},\check{\sigma}_{1})\in D_{1}$. 

\begin{example}
  Consider the example depicted in \reffig{fig:segment_4_rev}, where the 
  number of data points is $r=10$. 
  Let $k=5$ for example. 
  The partition $\{ D_{1},\dots,D_{5} \}$ of $D$ is given by 
  \begin{align*}
    (\check{\varepsilon}_{1},\check{\sigma}_{1}) ,\dots,
    (\check{\varepsilon}_{4},\check{\sigma}_{4}) 
    &\in D_{1} , \\
    (\check{\varepsilon}_{5},\check{\sigma}_{5}) ,\dots,
    (\check{\varepsilon}_{7},\check{\sigma}_{7}) 
    &\in D_{2} , \\
    (\check{\varepsilon}_{8},\check{\sigma}_{8}) ,\dots,
    (\check{\varepsilon}_{10},\check{\sigma}_{10}) 
    &\in D_{3} , 
  \end{align*}
  and $D_{4}=D_{5}=\emptyset$. 
  Correspondingly, 0-1 variables $t_{li}$ $(l=1,\dots,10;i=1,\dots,5)$ 
  take the values 
  \begin{alignat*}{5}
    & t_{1,1} = \dots = t_{4,1}  = 1 , 
    &{\quad}&
    t_{5,1} = \dots = t_{7,1} = 1 ,  
    &{\quad}&
    t_{8,1} = \dots = t_{10,1} = 1 ,  \\
    & t_{1,2} = \dots = t_{4,2}  = 0 , 
    &{\quad}&
    t_{5,2} = \dots = t_{7,2} = 1 , 
    &{\quad}&
    t_{8,2} = \dots = t_{10,2} = 1 ,  \\
    & t_{1,3} = \dots = t_{4,3}  = 0 , 
    &{\quad}&
    t_{5,3} = \dots = t_{7,3} = 0 , 
    &{\quad}&
    t_{8,3} = \dots = t_{10,3} = 1 ,  \\
    & t_{1,4} = \dots = t_{4,4}  = 0 , 
    &{\quad}&
    t_{5,4} = \dots = t_{7,4} = 0 , 
    &{\quad}&
    t_{8,4} = \dots = t_{10,4} = 0 ,  \\
    & t_{1,5} = \dots = t_{4,5}  = 0 , 
    &{\quad}&
    t_{5,5} = \dots = t_{7,5} = 0 ,
    &{\quad}&
    t_{8,5} = \dots = t_{10,5} = 0 . 
  \end{alignat*}
  Thus the partition of $D$ in \reffig{fig:segment_4_rev} 
  is expressed by using 0-1 variables. 
  It is worth noting that we have to give some relations 
  among $t_{l,1},\dots,t_{l,10}$ to ensure that, as explained in 
  section~\ref{sec:regression.problem}, 
  the data points belonging to the same $D_{i}$ should have consecutive 
  subscripts. 
  This motivates us to define the values of $t_{l,1},\dots,t_{l,10}$ 
  by \eqref{eq:def.binary.t.1}.  
  \finbox
\end{example}

It follows from \eqref{eq:def.binary.t.1} that, for each $l=1,\dots,r$, 
we have 
\begin{align}
  t_{l1} \ge t_{l2} \ge \dots \ge t_{lk} . 
  \label{eq:some.break.1}
\end{align}
Since the data points in $D_{i}$ should have consecutive subscripts, 
the constraints 
\begin{align*}
  t_{1i} \le t_{2i} \le \dots \le t_{ri} 
\end{align*}
should be satisfied for each $i=1,\dots,k$. 
Moreover, we can see that 
the number of nonempty ones among $D_{1},\dots,D_{k}$ (i.e., the number 
of straight lines used for data fitting) is 
equal to $t_{r1} + \dots + t_{rk}$.

The squared error of the regression model from data point 
$(\check{\varepsilon}_{l},\check{\sigma}_{l})$ can be expressed by using 
$t_{il}$ $(i=1,\dots,k)$ as follows. 
For each $l=1,\dots,r$, define $v_{l1},\dots,v_{lk} \in \Re$ by 
\begin{align}
  v_{li} = 
  \begin{dcases*}
    (\alpha_{i} \varepsilon_{l} + \beta_{i} \sigma_{l} - \gamma_{i})^{2} 
    & if $(t_{li},t_{l,i+1})=(1,0)$, \\
    0 
    & otherwise. 
  \end{dcases*}
  \label{eq:some.break.2}
\end{align}
Since \eqref{eq:def.binary.t.1} implies that 
$(t_{li},t_{l,i+1})=(1,0)$ if and only if 
$(\check{\varepsilon}_{l},\check{\sigma}_{l}) \in D_{i}$, 
we can write the squared error of $\ell_{i}$ from 
$(\check{\varepsilon}_{l},\check{\sigma}_{l})$ as 
$\sum_{i=1}^{k} v_{li}$. 
Therefore, the total error of the regression model is 
$\sum_{l=1}^{r} \sum_{i=1}^{k} v_{li}$. 
We next convert \eqref{eq:some.break.2} into some convex constraints. 
Observe that constraint \eqref{eq:some.break.1} implies that 
$(t_{li},t_{l,i+1})$ takes any one of $(1,0)$, $(0,0)$, and $(1,1)$. 
Letting $M$ be a sufficiently large constant, we see that 
\begin{align}
  M (1 - t_{li} + t_{l,i+1}) = 
  \begin{dcases*}
    0 & if $(t_{li},t_{l,i+1})=(1,0)$, \\
    M & if $(t_{li},t_{l,i+1})\in \{ (0,0),(1,1) \}$, \\
  \end{dcases*}
\end{align}
Therefore, $v_{li}$ satisfies \eqref{eq:some.break.2} if we 
minimize $v_{li}$ under the following constraints: 
\begin{align*}
  & v_{li} + M (1 - t_{li} + t_{l,i+1})
  \ge (\alpha_{i} \varepsilon_{l} + \beta_{i} \sigma_{l} - \gamma_{i})^{2} , \\
  & v_{li} \ge 0 . 
\end{align*}
It is worth noting that these two constraints are convex. 

As a consequence of the discussion above, 
the segmented least squares in section~\ref{sec:regression.problem} can 
be formulated as the following optimization problem: 
\begin{subequations}\label{P:regression}%
  \begin{alignat}{3}
    & \MIN  &{\quad}& 
    \sum_{l=1}^{r} \sum_{i=1}^{k} v_{li} + 
    \mu \sum_{i=1}^{k} t_{ri} \\
    & \ST && 
    t_{1i} \le \dots \le t_{ri}   , 
    \quad && i=1,\dots,k, \\
    & && 
    t_{l1} \ge \dots \ge t_{lk}   , 
    \quad && l=1,\dots,r, \\
    & && 
    t_{1,1} = 1  , 
    \label{P:regression.t11} \\
    & && 
    v_{li} + M (1 - t_{li} + t_{l,i+1}) \notag\\
    & && \qquad
    {}\ge (\alpha_{i} \varepsilon_{l} + \beta_{i} \sigma_{l} - \gamma_{i})^{2} , 
    \quad && l=1,\dots,r; \, i=1,\dots,k-1, \\
    & && 
    v_{lk} + M (1 - t_{lk}) \notag\\
    & && \qquad
    {}\ge (\alpha_{k} \varepsilon_{l} + \beta_{k} \sigma_{l} - \gamma_{k})^{2} , 
    \quad && l=1,\dots,r, \\
    & && 
    v_{l1},\dots,v_{lk} \ge 0 ,
    \quad && l=1,\dots,r, \\
    & && 
    t_{l1},\dots,t_{lk} \in \{ 0,1 \} , 
    \quad && l=1,\dots,r , \\
    & && 
    \beta_{i} = 1 , 
    \quad && i=1,\dots,k. 
    \label{P:regression.normalization}
  \end{alignat}
\end{subequations}
The optimization variables in this problem are 
$v_{li}$, $t_{li}$, $\alpha_{i}$, $\beta_{i}$, and $\gamma_{i}$ 
$(i=1,\dots,k; l=1,\dots,r)$. 
Since this optimization problem is a mixed-integer programming problem 
with convex constraints, we can find a global optimal solution with a 
branch-and-cut method. 
Particularly, since a convex quadratic constraint can be converted to a 
second-order cone constraint, 
we can recast this problem as an MISOCP problem, for which 
several well-developed solvers are available. 
Note that constraint \eqref{P:regression.normalization} is 
a normalization of parameters $\alpha_{i}$, $\beta_{i}$, 
and $\gamma_{i}$. 
Constraint \eqref{P:regression.t11} prevents a solution with 
$t_{li}=0$ and $v_{li}=0$ $(l=1,\dots,r;i=1,\dots,k)$ 
becomes optimal for problem~\eqref{P:regression}.\footnote{%
We have seen that \eqref{eq:constraint.for.D1} ensures 
$(\check{\varepsilon}_{1},\check{\sigma}_{1})\in D_{1}$. 
However, problem~\eqref{P:regression} does not involve 
$t_{1,2}=0$ as a constraint, because it is satisfied at the optimal 
solution by minimizing $\mu\sum_{i=1}^{k}t_{ri}$. }

\section{Expression and construction of uncertainty set}
\label{sec:uncertainty}

When we obtain the piecewise linear function in 
\reffig{fig:for_figure_prog3_segmented} by using the method developed in 
section~\ref{sec:regression}, we next construct the uncertainty set 
$C(\tau)$ shown in \reffig{fig:for_figure_prog3_seg_bound}. 
This section shows that 
condition $(\varepsilon,\sigma)\in C(\tau)$ is equivalently rewritten as 
constraints that can be handled within the framework of MILP. 
For ease of comprehension, section~\ref{sec:uncertainty.single} deals 
with the case that the result of the 
segmented least squares has a single breakpoint. 
We consider the general case in section~\ref{sec:uncertainty.many}. 
Section~\ref{sec:uncertainty.bisection} explains the procedure for 
determining the value of $\tau$.

\subsection{Single breakpoint case ($k=2$)}
\label{sec:uncertainty.single}

The lines obtained by the segmented least squares are denoted by 
\begin{align*}
  \ell_{i} : \quad
  \alpha_{i} \varepsilon + \beta_{i} \sigma = \gamma_{i} , 
  \quad i=1,2 . 
\end{align*}
Here, we introduce a normalization of parameters $\alpha_{i}$, 
$\beta_{i}$, and $\gamma_{i}$ as follows. 
Let $(\alpha_{i}^{*},\beta_{i}^{*},\gamma_{i}^{*})$ denote the optimal 
solution obtained by the method in section~\ref{sec:regression}. 
Then we normalize it as 
\begin{align}
  \alpha_{i} 
  = \frac{\alpha_{i}^{*}}{\| (\alpha_{i}^{*},\beta_{i}^{*}) \|} , 
  \quad
  \beta_{i} 
  = \frac{\beta_{i}^{*}}{\| (\alpha_{i}^{*},\beta_{i}^{*}) \|} , 
  \quad
  \gamma_{i} 
  = \frac{\gamma_{i}^{*}}{\| (\alpha_{i}^{*},\beta_{i}^{*}) \|} , 
  \label{eq:normalization.after.regression}
\end{align}
where $\| (\alpha_{i}^{*},\beta_{i}^{*}) \|$ denotes the Euclidean norm 
of vector $(\alpha_{i}^{*},\beta_{i}^{*})$. 
Moreover, without loss of generality we assume $\alpha_{i}<0$.

\begin{figure}[tbp]
  \centering
  \includegraphics[scale=0.70]{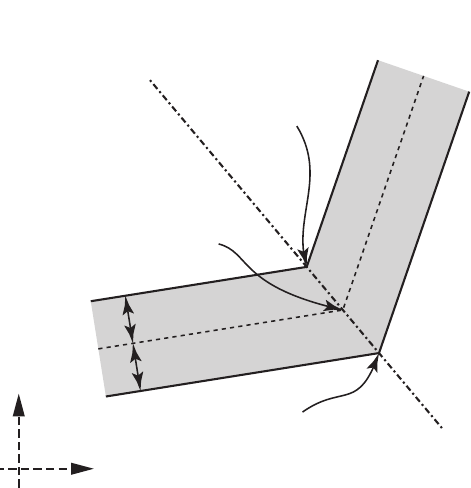}
  \begin{picture}(0,0)\small
    \put(-150,-50){
    \put(13,48){$\varepsilon$}
    \put(-15,76){$\sigma$}
    \put(11,111){$\ell_{1}^{+}$}
    \put(14,94){$\ell_{1}$}
    \put(15,77){$\ell_{1}^{-}$}
    \put(35,188){$\ell^{\rr{bd}}$}
    \put(115,198){$\ell_{2}^{+}$}
    \put(131,192){$\ell_{2}$}
    \put(147,186){$\ell_{2}^{-}$}
    \put(76,69){$(u^{-},v^{-})$}
    \put(34,130){$(u^{0},v^{0})$}
    \put(68,174){$(u^{+},v^{+})$}
    \put(38,89){$\tau$}
    \put(36,105){$\tau$}
    }
  \end{picture}
  \caption{Definition of set $C(\tau)$ (the gray colored part). }
  \label{fig:border_2}
\end{figure}
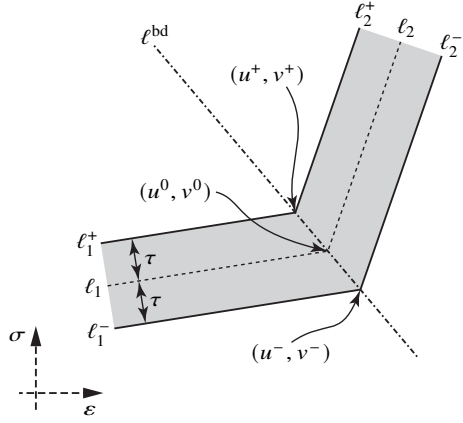

Define $C(\tau) \subset \Re^{2}$ as depicted in \reffig{fig:border_2}, 
where the border lines are obtained by shifting 
$\ell_{1}$ and $\ell_{2}$ by distance $\tau$; 
the precise definition of $C(\tau)$ will be given 
by \eqref{eq:single.set}. 
We use $\ell_{i}^{-}$ and $\ell_{i}^{+}$ to denote the lines consisting 
of the border of $C(\tau)$, where 
\begin{alignat}{3}
  \ell_{i}^{-} &: &{\quad}&
  \alpha_{i} \varepsilon + \beta_{i} \sigma = \gamma_{i} - \tau , 
  &{\quad}&  i=1,2 , 
  \label{eq:def.line.-} \\
  \ell_{i}^{+} &: &{\quad}&
  \alpha_{i} \varepsilon + \beta_{i} \sigma = \gamma_{i} + \tau ,
  &{\quad}&  i=1,2 . 
  \label{eq:def.line.+}
\end{alignat}

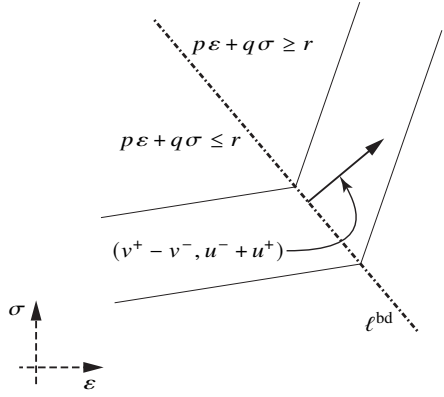
\begin{figure}[tbp]
  \centering
  \includegraphics[scale=0.70]{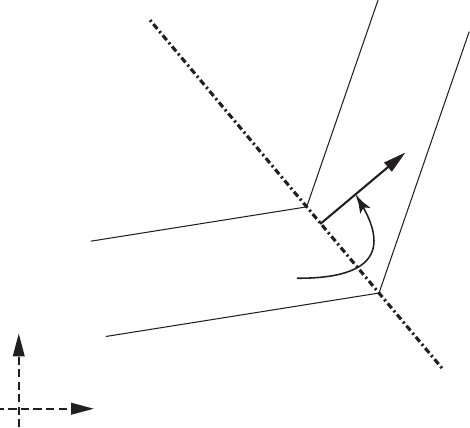}
  \begin{picture}(0,0)\small
    \put(-150,-50){
    \put(13,48){$\varepsilon$}
    \put(-15,76){$\sigma$}
    \put(27,140){$p\varepsilon + q\sigma \le r$}
    \put(54,176){$p\varepsilon + q\sigma \ge r$}
    \put(24,98){$(v^{+}-v^{-},u^{-}+u^{+})$}
    \put(119,70){$\ell^{\rr{bd}}$}
    }
  \end{picture}
  \caption{Case division in \eqref{eq:single.set}. }
  \label{fig:case_division}
\end{figure}

Let $(u^{-},v^{-})$ denote the intersection point of 
$\ell_{1}^{-}$ and $\ell_{2}^{-}$, and $(u^{+},v^{+})$ denote the 
intersection point of $\ell_{1}^{+}$ and $\ell_{2}^{+}$; 
see \reffig{fig:border_2}. 
We use $(u^{0},v^{0}) \in \Re^{2}$ to denote the intersection point of 
$\ell_{1}$ and $\ell_{2}$. 
Let $\ell^{\rr{bd}}$ denote the line passing through 
$(u^{-},v^{-})$, $(u^{0},v^{0})$, and $(u^{+},v^{+})$, 
and write it as 
\begin{align}
  \ell^{\rr{bd}} : \quad
  p \varepsilon + q \sigma = r . 
  \label{eq:single.p.q.r.1}
\end{align}
In the expression above, $p$, $q$, $r \in \Re$ are constants that will 
be computed concretely in the following. 
We see that $C(\tau)$ is a set of points 
$(\varepsilon,\sigma)\in\Re^{2}$ satisfying\footnote{%
Inclusion of condition $p \varepsilon + q \sigma = r$ in 
the both two cases of \eqref{eq:single.set} does not matter. } 
\begin{align}
  \begin{dcases*}
    |\alpha_{1} \varepsilon + \beta_{1} \sigma - \gamma_{1}| \le \tau 
    & if 
    $p \varepsilon + q \sigma \le r$, \\
    |\alpha_{2} \varepsilon + \beta_{2} \sigma - \gamma_{2}| \le \tau 
    & if 
    $p \varepsilon + q \sigma \ge r$; 
  \end{dcases*}
  \label{eq:single.set}
\end{align}
see \reffig{fig:border_2} and \reffig{fig:case_division}. 

We are now in position to compute the values of $p$, $q$, and $r$ 
for $\ell^{\rr{bd}}$ in \eqref{eq:single.p.q.r.1}; 
see \reffig{fig:case_division}. 
Observe that $\ell^{\rr{bd}}$ is in parallel with 
vector $(u^{+}-u^{-},v^{+}-v^{-})$ 
(namely, is orthogonal to vector $(v^{+}-v^{-},-u^{+}+u^{-})$), 
and passes through point $(u^{0},v^{0})$. 
Hence, it is written as 
\begin{align}
  (v^{+} - v^{-}) (\varepsilon - u^{0})
  - (u^{+} - u^{-}) (\sigma - v^{0}) = 0 . 
  \label{eq:single.p.q.r.5}
\end{align}
An elementary calculation shows that $u^{0}$, $v^{0}$, 
$u^{+} - u^{-}$, and $v^{+} - v^{-}$ are given by 
\begin{align}
  u^{0}
  &= \frac{\beta_{2}\gamma_{1} - \beta_{1}\gamma_{2}}
  {\alpha_{1}\beta_{2} - \alpha_{2}\beta_{1}} , 
  \label{eq:intersection.1} \\
  v^{0} 
  &= \frac{\alpha_{1}\gamma_{2} - \alpha_{2}\gamma_{1}}
  {\alpha_{1}\beta_{2} - \alpha_{2}\beta_{1}} , 
  \label{eq:intersection.2} \\
  u^{+} - u^{-} 
  &= \frac{2 (\beta_{2} - \beta_{1}) \tau}
  {\alpha_{1}\beta_{2} - \alpha_{2}\beta_{1}} , 
  \label{eq:u.diff.1} \\
  v^{+} - v^{-} 
  &= \frac{2 (\alpha_{1} - \alpha_{2}) \tau}
  {\alpha_{1}\beta_{2} - \alpha_{2}\beta_{1}} . 
  \label{eq:v.diff.1}
\end{align}
Substitute \eqref{eq:intersection.1}, \eqref{eq:intersection.2}, 
\eqref{eq:u.diff.1}, and \eqref{eq:v.diff.1} into 
\eqref{eq:single.p.q.r.5} to see that $p$, $q$, and $r$ 
in \eqref{eq:single.p.q.r.1} are obtained as 
\begin{align}
  p &= \alpha_{1} - \alpha_{2} , 
  \label{eq:single.p.q.r.2} \\
  q &= \beta_{1} - \beta_{2} , 
  \label{eq:single.p.q.r.3} \\
  r &= p u^{0} + q v^{0} . 
  \label{eq:single.p.q.r.4}
\end{align}

We next introduce a 0-1 variable $s$, and make it correspond to the case 
division in \eqref{eq:single.set} as\footnote{%
In a manner similar to \eqref{eq:single.set}, 
in \eqref{eq:single.division.2} 
condition $p \varepsilon + q \sigma = r$ corresponds to 
both $s=0$ and $s=1$, which is of no matter for the following formulations. }
\begin{alignat}{3}
  & s = 0 
  &{\quad}& \leftrightarrow &{\quad}& 
  p \varepsilon + q \sigma \le r , 
  \label{eq:single.division.1} \\
  & s = 1
  &{\quad}& \leftrightarrow &{\quad}& 
  p \varepsilon + q \sigma \ge r . 
  \label{eq:single.division.2}
\end{alignat}
Accordingly, the inequalities on the left side of \eqref{eq:single.set} 
correspond to $s$ as 
\begin{alignat}{3}
  & s = 0 
  &{\quad}& \leftrightarrow &{\quad}& 
  |\alpha_{1} \varepsilon + \beta_{1} \sigma - \gamma_{1}| \le \tau , 
  \label{eq:single.division.3} \\
  & s = 1
  &{\quad}& \leftrightarrow &{\quad}& 
  |\alpha_{2} \varepsilon + \beta_{2} \sigma - \gamma_{2}| \le \tau . 
  \label{eq:single.division.4}
\end{alignat}
Thus, \eqref{eq:single.set} is equivalently restated as 
\eqref{eq:single.division.1}, \eqref{eq:single.division.2}, 
\eqref{eq:single.division.3}, and \eqref{eq:single.division.4}. 
The remaining concern is to reduce these four conditions to some linear 
inequalities. 
To this end, we use a sufficiently large constant $M$. 
Then conditions \eqref{eq:single.division.1} and 
\eqref{eq:single.division.2} are equivalently rewritten as 
\begin{align}
  p \varepsilon + q \sigma 
  &\le r + M s , 
  \label{eq:single.break.1}\\
  p \varepsilon + q \sigma 
  &\ge r - M (1-s) . 
  \label{eq:single.break.2}
\end{align}
We also see that conditions \eqref{eq:single.division.3} and 
\eqref{eq:single.division.4} are equivalently rewritten as 
\begin{align}
  |\alpha_{1} \varepsilon + \beta_{1} \sigma - \gamma_{1}| 
  &\le \tau + M s , 
  \label{eq:single.break.3} \\
  |\alpha_{2} \varepsilon + \beta_{2} \sigma - \gamma_{2}| 
  &\le \tau + M (1-s) . 
  \label{eq:single.break.4}
\end{align}

As a consequence of this section, we obtain 
\begin{align*}
  (\varepsilon,\sigma) \in C(\tau) 
  \quad\Leftrightarrow\quad
  \exists s \in \{ 0, 1 \}  : 
  \eqref{eq:single.break.1}, \eqref{eq:single.break.2}, 
  \eqref{eq:single.break.3}, \eqref{eq:single.break.4} . 
\end{align*}
It is worth noting that the constraints on the right side of this 
expression can be treated within the framework of MILP.

\subsection{General case $(k \ge 3)$}
\label{sec:uncertainty.many}

This section generalizes the analysis in 
section~\ref{sec:uncertainty.single} to 
the case in which the result of the 
segmented least squares has $k$ lines with $k-1$ breakpoints. 

As the output of the segmented least squares, we obtain 
lines $\ell_{1},\dots,\ell_{k}$. 
We use the normalization of 
$\alpha_{i}$, $\beta_{i}$, $\gamma_{i}$ $(i=1,\dots,k)$ in 
\eqref{eq:normalization.after.regression}, and assume 
$\alpha_{i} < 0$ without loss of generality. 
Let $(u^{0}_{i},v^{0}_{i})$ denote the intersection point 
of $\ell_{i}$ and $\ell_{i+1}$ $(i=1,\dots,k-1)$. 
A simple calculation shows 
\begin{align}
  (u_{i}^{0},v_{i}^{0}) 
  = \Bigl(
  \frac{\beta_{i+1}\gamma_{i} - \beta_{i}\gamma_{i+1}}
  {\alpha_{i}\beta_{i+1} - \alpha_{i+1}\beta_{i}} , 
  \frac{\gamma_{i+1}\alpha_{i} - \gamma_{i}\alpha_{i+1}}
  {\alpha_{i}\beta_{i+1} - \alpha_{i+1}\beta_{i}} 
  \Bigr) . 
  \label{eq:u.v.i.0}
\end{align}
The lines consisting of the border of $C(\tau)$, denoted by 
$\ell_{i}^{-}$ and $\ell_{i}^{+}$ $(i=1,\dots,k)$, are defined in the 
same manner as \eqref{eq:def.line.-} and \eqref{eq:def.line.+}. 
Define $p_{i}$, $q_{i}$, $r_{i} \in \Re$ $(i=1,\dots,k-1)$ by 
\begin{align}
  p_{i} &= \alpha_{i} - \alpha_{i+1} , 
  \label{eq:many.lines.p} \\
  q_{i} &= \beta_{i} - \beta_{i+1} , 
  \label{eq:many.lines.q} \\
  r_{i} &= p_{i} u^{0}_{i} + q_{i} v^{0}_{i} . 
  \label{eq:many.lines.r}
\end{align}
Lines defined by $p_{i} \varepsilon + q_{i} \sigma = r_{i}$ 
$(i=1,\dots,k-1)$ play a role 
similar to $\ell^{\rr{bd}}$ in section~\ref{sec:uncertainty.single}; 
see \reffig{fig:border_2} and \reffig{fig:case_division}. 
For simplicity of notation, let 
\begin{alignat}{3}
  p_{0} &= 0, 
  &\quad
  q_{0} &= 0 , 
  &\quad
  r_{0} &= 0 , 
  \label{eq:many.lines.0} \\
  p_{k} &= 0, 
  &\quad
  q_{k} &= 0 , 
  &\quad
  r_{k} &= 0 . 
  \label{eq:many.lines.k}
\end{alignat}
We see that $C(\tau)$ is the set of points satisfying 
\begin{align}
  |\alpha_{i} \varepsilon + \beta_{i} \sigma - \gamma_{i}| \le \tau 
  \quad \Leftarrow \quad
  & p_{i-1} \varepsilon + q_{i-1} \sigma \ge r_{i-1}, \
  p_{i} \varepsilon + q_{i} \sigma \le r_{i} ,  \notag\\
  & \qquad\qquad\qquad i=1,\dots,k . 
  \label{eq:multiple.set}
\end{align}

We next introduce $0$-$1$ variables 
$s_{1},\dots,s_{k-1} \in \{ 0,1 \}$ to express the case divisions in 
\eqref{eq:multiple.set} as 
\begin{align}
  & s_{1} = \dots = s_{i-1}= 1 , 
  s_{i} = \dots = s_{k-1} = 0 \notag \\
  & \qquad \leftrightarrow\quad
  p_{i-1} \varepsilon + q_{i-1} \sigma \ge r_{i-1}, \
  p_{i} \varepsilon + q_{i} \sigma \le r_{i} 
  \label{eq:multiple.correspondence.1}
\end{align}
for each $i=1,\dots,k-1$. 
In this expression, we see that the $0$-$1$ variables satisfy 
\begin{align}
  s_{1} \ge s_{2} \ge \dots \ge s_{k-1} . 
  \label{eq:multiple.constraint.0}
\end{align}
The inequality on the left side of \eqref{eq:multiple.set} is linked 
to the 0-1 variables as 
\begin{align}
  & s_{1} = \dots = s_{i-1}= 1 , 
  s_{i} = \dots = s_{k-1} = 0 \notag \\
  & \qquad \leftrightarrow\quad
  |\alpha_{i} \varepsilon + \beta_{i} \sigma - \gamma_{i}| 
  \le \tau
  \label{eq:multiple.correspondence.2}
\end{align}
for each $i=1,\dots,k$. 
Thus, $C(\tau)$ in \eqref{eq:multiple.set} is equivalently expressed 
as \eqref{eq:multiple.correspondence.1} and \eqref{eq:multiple.correspondence.2}. 
Observe that the relation in \eqref{eq:multiple.correspondence.1} can be 
rewritten as the linear inequalities 
\begin{align}
  p_{i-1} \varepsilon + q_{i-1} \sigma 
  & \ge r_{i-1} 
  - M (i-1 - s_{1} - \dots - s_{i-1} + s_{i} + \dots + s_{k-1}) , 
  \notag\\
  & \qquad\qquad i=1,\dots,k , 
  \label{eq:multiple.constraint.1} \\
  p_{i} \varepsilon + q_{i} \sigma 
  & \le r_{i} 
  + M (i-1 - s_{1} - \dots - s_{i-1} + s_{i} + \dots + s_{k-1}) , 
  \notag\\
  & \qquad\qquad i=1,\dots,k . 
  \label{eq:multiple.constraint.2}
\end{align}
Similarly, the relation in \eqref{eq:multiple.correspondence.2} can be 
rewritten as the following linear inequalities: 
\begin{align}
  |\alpha_{i} \varepsilon + \beta_{i} \sigma - \gamma_{i}| 
  & \le \tau + M (i-1 - s_{1} -\dots- s_{i-1} + s_{i} + \dots + s_{k-1}) , 
  \notag\\
  & \qquad\qquad i=1,\dots,k . 
  \label{eq:multiple.constraint.7}
\end{align}

As a consequence, point $(\varepsilon,\sigma) \in \Re^{2}$ satisfies 
\eqref{eq:multiple.set} (i.e., 
$(\varepsilon,\sigma) \in C(\tau)$) if and only if there exist 
$(s_{1},\dots,s_{k-1}) \in \{ 0,1 \}^{k}$ satisfying 
\eqref{eq:multiple.constraint.0}, \eqref{eq:multiple.constraint.1}, 
\eqref{eq:multiple.constraint.2}, and \eqref{eq:multiple.constraint.7}.

\subsection{Determination of $\tau$}
\label{sec:uncertainty.bisection}

This section describes the procedure for determining 
the value of $\tau$ for $C(\tau)$. 
Recall that the data set $D$ consists of $r$ data points. 
When we choose the values of $\epsilon$ and $\delta$ in 
\eqref{eq:confidence.inequality.1}, 
$\tilde{p}$ is determined as the smallest integer satisfying 
\eqref{eq:def.tilde.p}. 
Then $\tau$ is to be determined so that the number of data points that 
are included in $C(\tau)$ is equal to $\tilde{p}$. 
We can use a bisection method for finding this value of $\tau$. 

The bisection method demands the following preparation. 
We have the output, 
$(\alpha_{i},\beta_{i},\gamma_{i})$ $(i=1,\dots,k)$, of the 
segmented least squares. 
Compute $u_{i}^{0}$, $v_{i}^{0}$ $(i=1,\dots,k-1)$ 
by \eqref{eq:u.v.i.0}. 
Moreover, compute $p_{i}$, $q_{i}$, $r_{i}$ $(i=0,1,2,\dots,k)$ 
by \eqref{eq:many.lines.p}, \eqref{eq:many.lines.q}, 
\eqref{eq:many.lines.r}, \eqref{eq:many.lines.0}, 
and \eqref{eq:many.lines.k}. 
Define $D_{i} \subset \Re^{2}$ $(i=1,\dots,k)$ by 
\begin{align*}
  D_{i} 
  = \{ (\varepsilon,\sigma) \in \Re^{2}
  \mid
  p_{i-1} \varepsilon + q_{i-1} \sigma \ge r_{i-1} , \
  p_{i} \varepsilon + q_{i} \sigma \le r_{i} \} . 
\end{align*}
Then, for each $l=1,\dots,k$, 
find $i_{l} \in \{ 1,\dots,k \}$ such that 
$(\check{\varepsilon}_{l},\check{\sigma}_{l}) \in D_{i_{l}}$ holds. 

\refalg{alg:bisection.1} describes the bisection method for 
computing $\tau$. 

\begin{algorithm}
  \caption{Bisection method for computing $\tau$. }
  \label{alg:bisection.1}
  \begin{algorithmic}[1]
    \Require
    $i_{l}$ $(l=1,\dots,r)$, $\tilde{p}$, 
    $\tau_{\rr{max}}>0$, $\epsilon_{\rr{bi}}>0$. 
    \State
    $\tau \gets \tau_{\rr{max}}$, 
    $\tau_{\rr{min}} \gets 0$. 
    \While{$\tau_{\rr{max}} - \tau_{\rr{min}} \ge \epsilon_{\rr{bi}}$}
    \State
    $\tau \gets (\tau_{\rr{max}} + \tau_{\rr{min}})/2$, 
    $p \gets 0$.
    \For{$l=1,\dots,r$}
    \If{$|\alpha_{i_{l}} \check{\varepsilon}_{l} 
    + \beta_{i_{l}} \check{\sigma}_{l} - \gamma_{i_{l}}| 
    \le \tau$}
    \State
    $p \gets p+1$. 
    \EndIf
    \EndFor
    \If{$p \ge \tilde{p}$}
    \State
    $\tau_{\rr{max}} \gets \tau$. 
    \Else
    \State
    $\tau_{\rr{min}} \gets \tau$. 
    \EndIf
    \EndWhile
  \end{algorithmic}
\end{algorithm}

\section{Formulation for trusses}
\label{sec:truss}

By adopting $C(\tau)$ formulated in section~\ref{sec:uncertainty}, 
we can compute lower and upper bounds for the structural response that 
satisfy \eqref{eq:confidence.inequality.1}. 
In this section, taking trusses for example, we show that the optimization 
problems for computing these lower and upper bounds can be recast as  
MILP problems. 
It is worth noting that we can solve an MILP problem globally with, e.g., 
a branch-and-cut method. 
Guaranteeing the global optimality is crucial to finding a bound for the 
structural response, because a local optimal solution in general 
underestimates the structural response; 
the global optimality ensures \eqref{eq:confidence.inequality.1}. 

Consider a truss that undergoes small deformation. 
The constitutive law relates the uniaxial strain 
$\varepsilon\in\Re$ to the uniaxial stress $\sigma\in \Re$. 
Let $d$ denote the number of degrees of freedom of the nodal 
displacements of the truss. 
We use $\bi{u} \in \Re^{d}$ and $\bi{f} \in \Re^{d}$ to denote the 
vectors of the nodal displacements and the nodal external forces, 
respectively. 
Let $m$ denote the number of truss members. 
For member $e$ $(e=1,\dots,m)$, we use 
$\sigma_{e} \in \Re$ and $\varepsilon_{e} \in \Re$ to denote its stress 
and strain, respectively, and write 
$\bi{\varepsilon} = (\varepsilon_{e}) \in \Re^{m}$ and 
$\bi{\sigma} = (\sigma_{e}) \in \Re^{m}$. 
The compatibility relations and the force-balance equations can be 
written in the forms 
\begin{align*}
  \bi{\varepsilon}  &= L \bi{u} , \\
  N \bi{\sigma}  &= \bi{f} , 
\end{align*}
where $L \in \Re^{m\times d}$ and $N \in \Re^{d \times m}$ are 
constant matrices. 

We obtain the uncertainty set, $C(\tau)$, by using the method 
presented in section~\ref{sec:uncertainty}. 
Then we can see from the fact reviewed in 
section~\ref{sec:framework.bound} that $C(\tau)$ satisfies 
\begin{align*}
  \rr{P}_{F}\bigl\{ 
  \rr{P}\{ (\varepsilon,\sigma) \in C(\tau) \} \ge 1-\epsilon 
  \bigr\} \ge 1-\delta . 
\end{align*}
Let $q(\bi{u},\bi{\sigma})$ denote the quantity of interest. 
Then a lower bound for the quantity of interest, $\underline{q}$, can be 
obtained as the optimal value of the following optimization problem: 
\begin{subequations}\label{P:truss.lower}%
  \begin{alignat}{3}
    & \MIN
    &{\quad}& 
    q(\bi{u},\bi{\sigma})  \\
    & \ST && 
    \bi{\varepsilon} = L \bi{u} , \\
    & &&
    N \bi{\sigma} = \bi{p} , \\
    & &&
    (\varepsilon_{e}, \sigma_{e}) \in C(\tau) , 
    \quad e=1,\dots,m . 
    \label{eq:uncertainty.inclusion.2}
  \end{alignat}
\end{subequations}
It is crucial in the formulation above that constraint 
\begin{align}
  (\varepsilon_{e}, \sigma_{e}) \in C(\tau) 
  \label{eq:uncertainty.inclusion.1}
\end{align}
for each $e=1,\dots,m$ can be equivalently rewritten into 
some linear inequalities by introducing some $0$-$1$ variables. 

We next present the concrete reformulation of 
constraint \eqref{eq:uncertainty.inclusion.1}. 
As an example, suppose that the number of straight lines of the output 
of the segmented least squares is $k=3$ (i.e., the number of breakpoints 
is two). 
We use the result established in section~\ref{sec:uncertainty.many}. 
For each member $e$ $(e=1,\dots,m)$, we introduce two $0$-$1$ variables, 
$s_{e1}$ and $s_{e2}$. 
Then we can reformulate \eqref{eq:uncertainty.inclusion.1} equivalently 
as follows: 
\begin{alignat}{3}
  & s_{e1} \ge s_{e2} , 
  \label{eq:truss.explicit.1} \\
  & p_{1} \varepsilon_{e} + q_{1} \sigma_{e} 
  \le r_{1}  + M (s_{e1} + s_{e2}) , 
  \label{eq:truss.explicit.2} \\
  & p_{1} \varepsilon_{e} + q_{1} \sigma_{e} 
  \ge r_{1} - M (1 - s_{e1} + s_{e2}) , 
  \label{eq:truss.explicit.3} \\
  & p_{2} \varepsilon_{e} + q_{2} \sigma_{e} 
  \le r_{2}  + M (1 - s_{e1} + s_{e2}) , 
  \label{eq:truss.explicit.4} \\
  & p_{2} \varepsilon_{e} + q_{2} \sigma_{e} 
  \ge r_{2} - M (2 - s_{e1} - s_{e2}) , 
  \label{eq:truss.explicit.5} \\
  & |\alpha_{1} \varepsilon_{e} + \beta_{1} \sigma_{e} - \gamma_{1}| 
  \le d + M (s_{e1}+s_{e2}) , 
  \label{eq:truss.explicit.6} \\
  & |\alpha_{2} \varepsilon_{e} + \beta_{2} \sigma_{e} - \gamma_{2}| 
  \le d + M (1-s_{e1}+s_{e2}) , 
  \label{eq:truss.explicit.7} \\
  & |\alpha_{3} \varepsilon_{e} + \beta_{3} \sigma_{e} - \gamma_{3}| 
  \le d + M (2-s_{e1}-s_{e2}) , 
  \label{eq:truss.explicit.8} \\
  & s_{e1}, s_{e2} \in \{ 0,1 \} . 
  \label{eq:truss.explicit.9}
\end{alignat}
Thus, we can replace constraint \eqref{eq:uncertainty.inclusion.2} of 
problem \eqref{P:truss.lower} with 
\eqref{eq:truss.explicit.1}, \eqref{eq:truss.explicit.2}, 
\eqref{eq:truss.explicit.3}, \eqref{eq:truss.explicit.4}, 
\eqref{eq:truss.explicit.5}, \eqref{eq:truss.explicit.6}, 
\eqref{eq:truss.explicit.7}, \eqref{eq:truss.explicit.8}, and 
\eqref{eq:truss.explicit.9} for each $e=1,\dots,m$. 
The optimization problem is thence an MILP problem when the objective 
function $q$ is a linear function (that is often the case as seen in 
section~\ref{sec:ex}). 
Accordingly, we can find the global optimal solution of 
problem \eqref{P:truss.lower}, which corresponds to a conservative 
prediction $\underline{q}$ satisfying \eqref{eq:confidence.inequality.1}. 

An upper bound for the quantity of interest, $\overline{q}$, can be 
found by maximizing $q$ under the same constraints.

\section{Numerical examples}
\label{sec:ex}

This section demonstrates three numerical examples.\footnote{%
The main MATLAB codes and the data sets used in the analysis of this 
section are available at 
\url{https://github.com/ykanno22/rel_comp_segmented/}. } 
The proposed method was implemented on MATLAB ver.~23.2. 
We used CPLEX ver.~12.9~\citep{cplex} to solve MILP and MISOCP problems. 
For MISOCP problems, we set the \url{MIQCP strategy} parameter of CPLEX 
so that linear programming relaxations are adopted, rather than 
relaxations with convex quadratic constraints. 
The numerical experiments were conducted on 
a 2.6{\,}GHz Intel Core i7 processor with 32{\,}GB RAM.

\subsection{Truss example}
\label{sec:ex.truss}

\begin{figure}[tbp]
  \centering
  \includegraphics[scale=1.00]{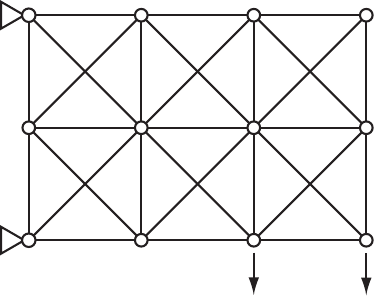}
  \begin{picture}(0,0)
    \put(-148,-50){
    \put(3,177){{\footnotesize $1$}}
    \put(57,177){{\footnotesize $2$}}
    \put(112,177){{\footnotesize $3$}}
    \put(3,122){{\footnotesize $4$}}
    \put(57,122){{\footnotesize $5$}}
    \put(112,122){{\footnotesize $6$}}
    \put(6,68){{\footnotesize $7$}}
    \put(60,68){{\footnotesize $8$}}
    \put(114,68){{\footnotesize $9$}}
    \put(-30,154){{\footnotesize $10$}}
    \put(24,154){{\footnotesize $11$}}
    \put(78,154){{\footnotesize $12$}}
    \put(144,154){{\footnotesize $13$}}
    \put(-30,100){{\footnotesize $14$}}
    \put(24,100){{\footnotesize $15$}}
    \put(78,100){{\footnotesize $16$}}
    \put(144,100){{\footnotesize $17$}}
    \put(-12,161){{\footnotesize $18$}}
    \put(42,161){{\footnotesize $19$}}
    \put(96,161){{\footnotesize $20$}}
    \put(-12,107){{\footnotesize $21$}}
    \put(42,107){{\footnotesize $22$}}
    \put(96,107){{\footnotesize $23$}}
    \put(-4,139){{\footnotesize $24$}}
    \put(50,139){{\footnotesize $25$}}
    \put(104,139){{\footnotesize $26$}}
    \put(-4,85){{\footnotesize $27$}}
    \put(50,85){{\footnotesize $28$}}
    \put(104,85){{\footnotesize $29$}}
    }
  \end{picture}
  \caption{Problem setting of the truss example.}
  \label{fig:x3_y2_truss}
\end{figure}

\begin{figure}[tbp]
  \centering
  \includegraphics[scale=0.50]{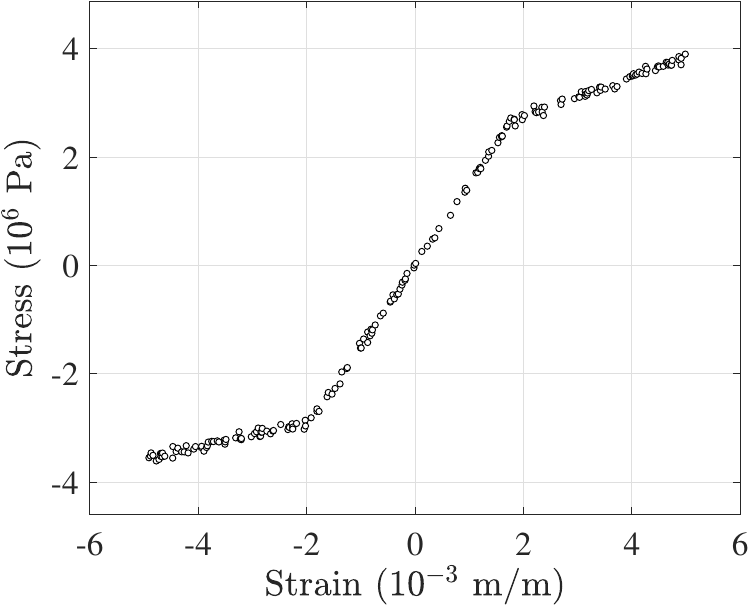}
  \caption{Data set for the truss example. }
  \label{fig:prog3_tri-modulus_1d_data_set}
\end{figure}

\begin{figure}[tbp]
  \centering
  \includegraphics[scale=0.50]{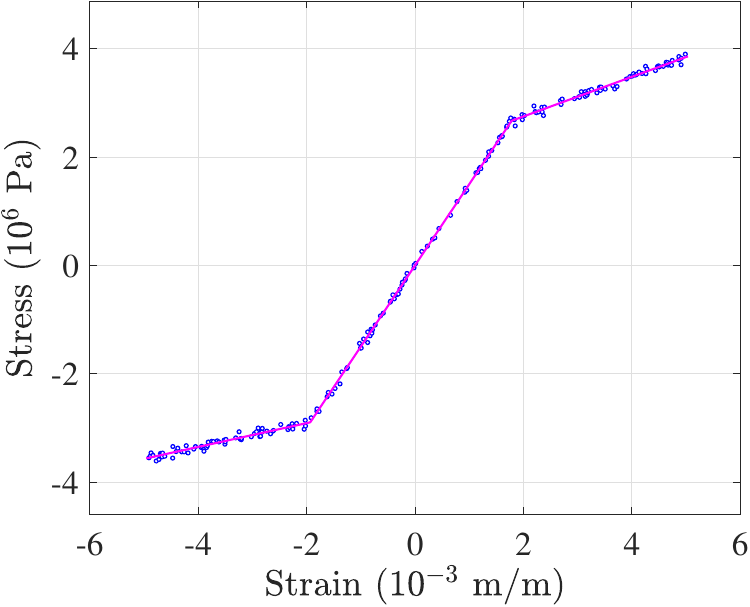}
  \caption{Result of the segmented least squares for the data set in 
  \reffig{fig:prog3_tri-modulus_1d_data_set}. }
  \label{fig:prog3_segmented_regress_result}
\end{figure}

\begin{figure}[tbp]
  \centering
  \includegraphics[scale=0.50]{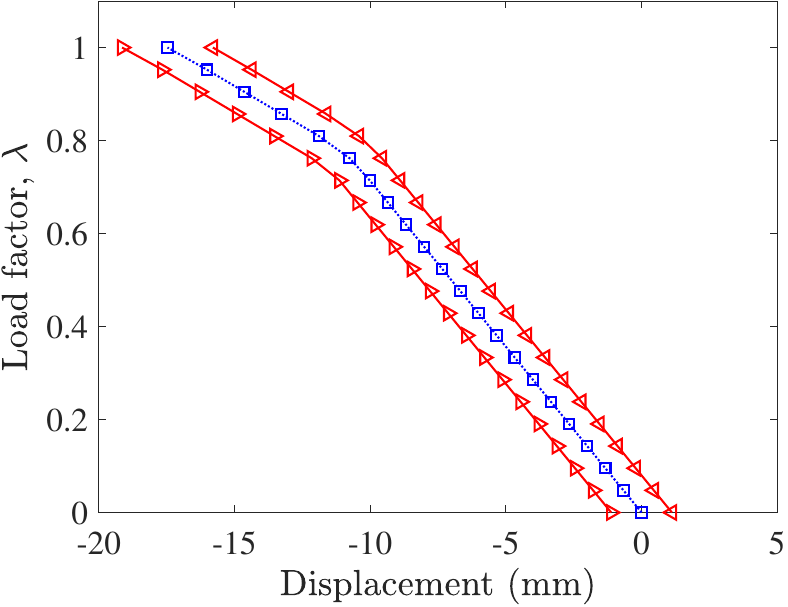}
  \caption{Bounds for the nodal displacement of the truss. 
  ``$\triangleleft$'' and ``$\triangleright$'' denote the upper and lower 
  bounds obtained by the proposed method, respectively; 
  ``$\square$'' denotes the reference solution. }
  \label{fig:prog3_displ_analysis_loop}
\end{figure}

\begin{figure}[tbp]
  \centering
  \includegraphics[scale=0.50]{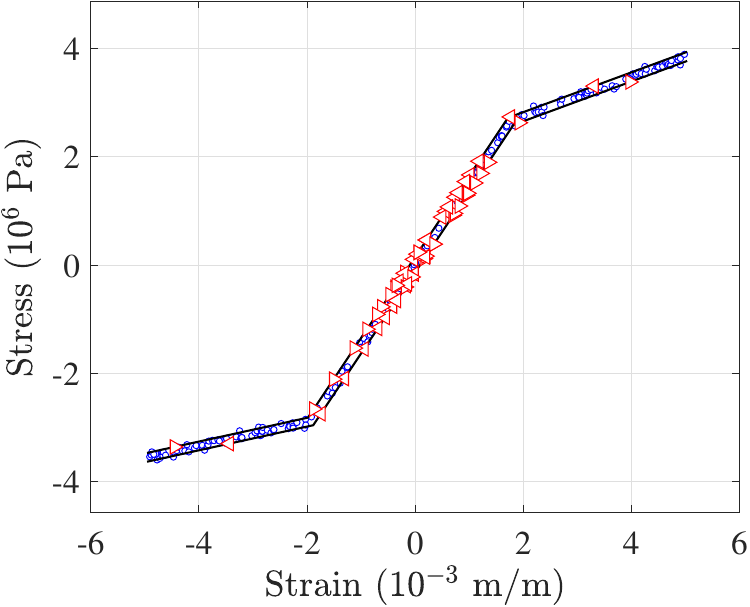}
  \caption{Member stresses--strains corresponding to the obtained 
  solutions at $\lambda=1$ in \reffig{fig:prog3_displ_analysis_loop}. 
  Boundary of the uncertainty set $C(\tau)$ is also depicted. 
  ``$\triangleleft$'' and ``$\triangleright$'' denote the member 
  stresses corresponding to the upper and lower bound solutions, 
  respectively. }
  \label{fig:prog3_displ_analysis_result}
\end{figure}

\begin{figure}[tbp]
  \centering
  \includegraphics[scale=0.50]{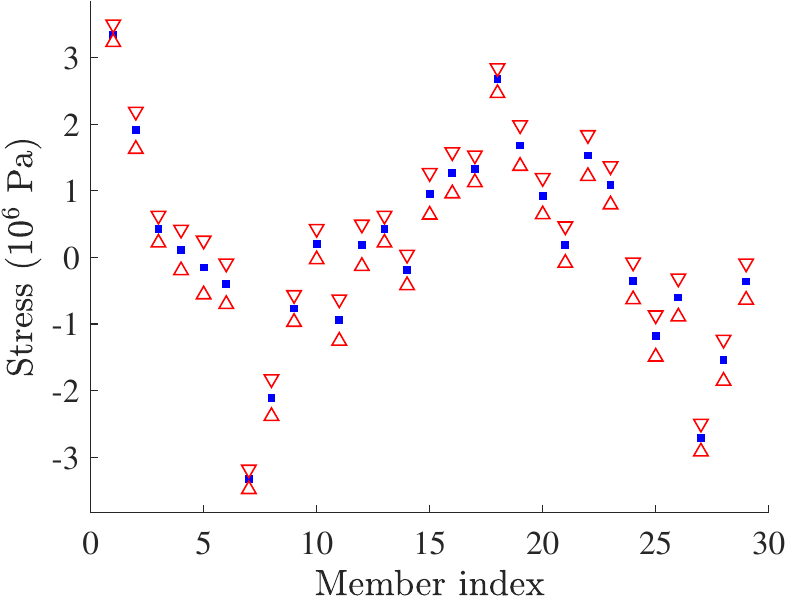}
  \caption{Bounds for the member stresses of the truss example. 
  ``$\triangledown$'' and ``$\vartriangle$'' denote the upper and lower 
  bounds obtained by the proposed method, respectively. 
  ``{\footnotesize $\blacksquare$}'' denotes the reference solution. }
  \label{fig:prog3_stress_all_bound}
\end{figure}

\begin{figure}[tbp]
  \centering
  \subfloat[]{
  \label{fig:prog3_stress_delta_variation}
  \includegraphics[scale=0.50]{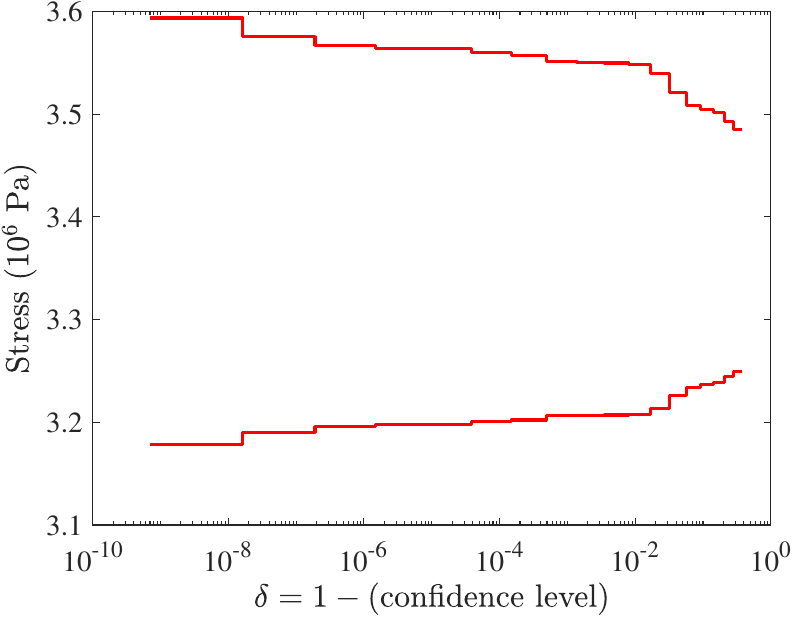}
  }
  \hfill
  \subfloat[]{
  \label{fig:prog3_stress_epsilon_variation}
  \includegraphics[scale=0.50]{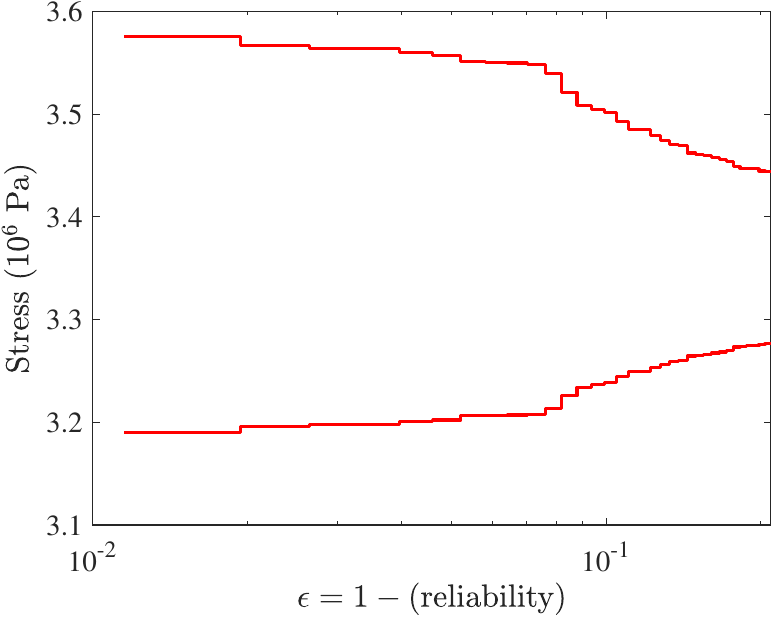}
  }
  \caption[]{Variations of the upper and lower bounds for the stress of 
  member 1 in the truss example. 
  \subref{fig:prog3_stress_delta_variation}~$\epsilon=0.1$; and 
  \subref{fig:prog3_stress_epsilon_variation}~$\delta=0.1$. }
  \label{fig:prog3_stress_variation}
\end{figure}

Consider a planar truss shown in \reffig{fig:x3_y2_truss}. 
This truss consists of $m=29$ members and has $d=20$ degrees of freedom 
of the nodal displacements. 
The lengths of horizontal and vertical members are $1\,\mathrm{m}$. 
The cross-sectional area of each member is $1000\,\mathrm{mm^{2}}$. 
We apply the external vertical downward forces of 
$2.1\lambda\,\mathrm{kN}$ at the bottom two nodes as shown in 
\reffig{fig:x3_y2_truss}, where $\lambda\in\Re$ is the load factor. 

Suppose that we are given the material data set shown in 
\reffig{fig:prog3_tri-modulus_1d_data_set}, which consists of 
$r=200$ data points. 
\reffig{fig:prog3_segmented_regress_result} shows the optimal solution 
of the segmented least squares with $k=5$, $\mu=2.0$, and 
$M=1.0\times 10^{3}$. 
We can observe that the penalty for the number of lines works 
substantially: Although five straight lines can be used for fitting the 
data, the optimal solution uses only three lines and fits the data very well. 

In accordance with the result of the segmented least squares 
in \reffig{fig:prog3_segmented_regress_result}, 
we put $k=3$ in the following. 
We specify the target reliability and confidence level 
as $1-\epsilon=0.9$ and $1-\delta=0.9$, respectively. 
The uncertainty set, $C(\tau)$, is then determined by 
\refalg{alg:bisection.1} with $\epsilon_{\rr{bi}}=10^{-7}$. 
To find a bound for the structural response, 
we employ the MILP formulation in section~\ref{sec:truss} 
with $M=1.0\times 10^{2}$. 
As for the quantity of interest, we first adopt the vertical 
displacement of the rightmost bottom node of the truss. 
\reffig{fig:prog3_displ_analysis_loop} reports the solutions obtained by 
the proposed method, where the obtained bounds are shown for several 
different values of the load factor, $\lambda$. 
As for a reference solution, we use the 
result of the segmented least squares 
in \reffig{fig:prog3_segmented_regress_result} as the constitutive law, 
and perform the conventional equilibrium analysis.\footnote{%
We adopt MATLAB built-in function \texttt{fsolve} as a solver for 
a system of nonlinear equations. } 
These reference solutions are also shown in 
\reffig{fig:prog3_displ_analysis_loop} as ``$\square$''. 
For each value of $\lambda$, we can observe that the reference solution 
belongs to the obtained interval, 
which supports the rationality of the proposed method. 
\reffig{fig:prog3_displ_analysis_result} shows the boundary of $C(\tau)$, 
and plots $(\varepsilon_{e},\sigma_{e})$ $(e=1,\dots,m)$ corresponding 
to the obtained upper and lower bounds for $\lambda=1$. 
We can confirm $(\varepsilon_{e},\sigma_{e}) \in C(\tau)$ for all 
the members. 

We next focus attention to the member stresses at $\lambda=1$. 
\reffig{fig:prog3_stress_all_bound} shows the upper and lower 
bonds obtained by the proposed method as well as the reference 
solutions, where the member indices are 
defined in \reffig{fig:x3_y2_truss}. 

\reffig{fig:prog3_stress_delta_variation} reports the variation of the 
interval bound for the stress of member~1 with respect to the confidence 
level $1-\delta$, with $\epsilon=0.1$ being fixed. 
As the required confidence level becomes higher, the interval 
guaranteeing the reliability $1-\epsilon$ becomes wider. 
Similarly, \reffig{fig:prog3_stress_epsilon_variation} reports the 
variation of the interval bound with respect to the target reliability 
$1-\epsilon$, with $\delta$ being fixed. 
As the target reliability becomes higher, 
the interval with the required confidence level becomes wider.

\subsection{Frame example}
\label{sec:ex.frame}

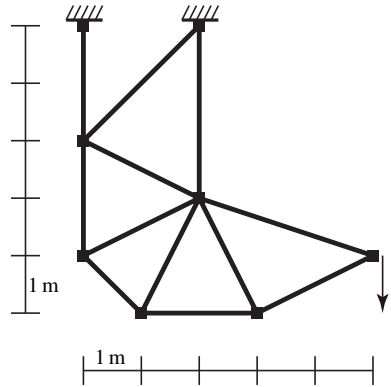
\begin{figure}[tbp]
  \centering
  \includegraphics[scale=0.60]{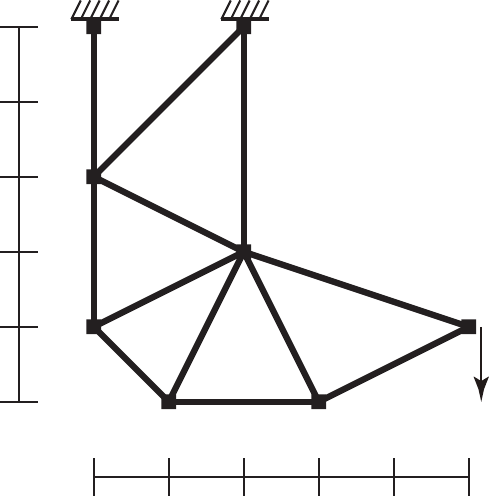}
  \begin{picture}(0,0)
    \put(-150,-50){
    \put(38,58){{\footnotesize $1\,\mathrm{m}$}}
    \put(13,85){{\footnotesize $1\,\mathrm{m}$}}
    }
  \end{picture}
  \caption{Problem setting of the frame example.}
  \label{fig:prog4_frame}
\end{figure}

\begin{figure}[tbp]
  \centering
  \includegraphics[scale=0.50]{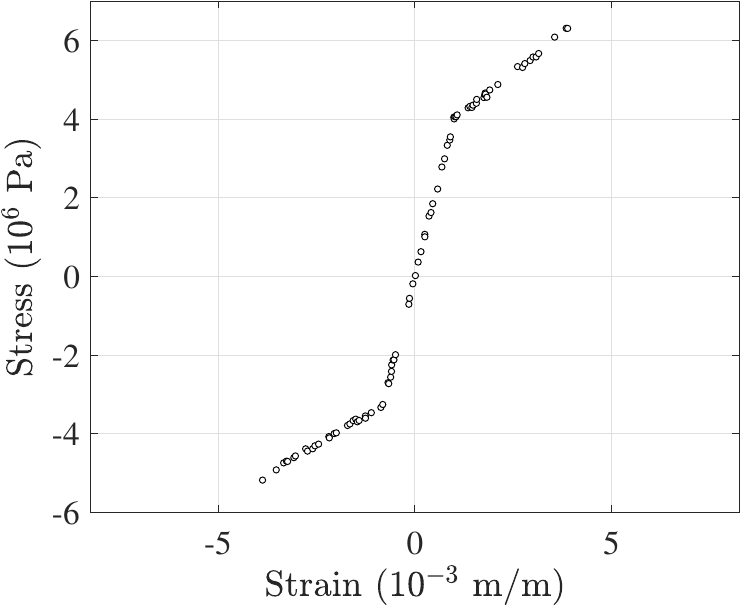}
  \caption{Data set for the frame example.}
  \label{fig:prog4_tri-modulus_1d_data_set}
\end{figure}

\begin{figure}[tbp]
  \centering
  \includegraphics[scale=0.50]{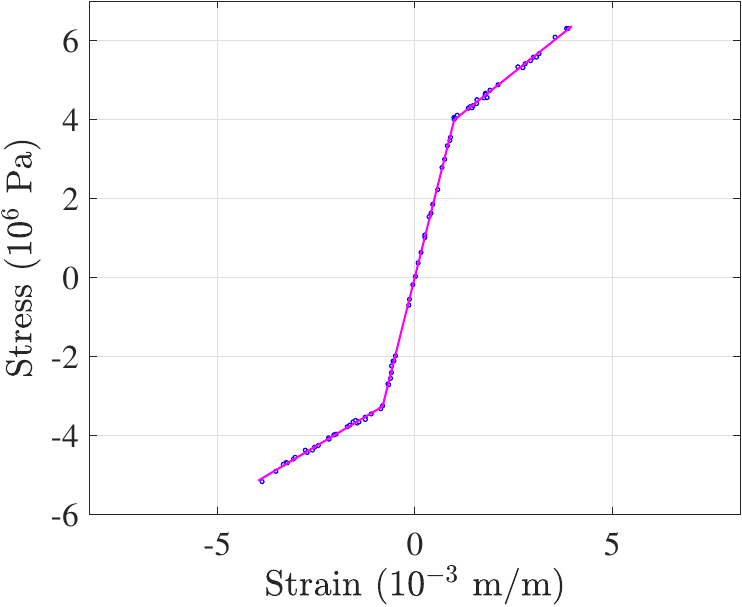}
  \caption{Result of the segmented least squares for the data set in 
  \reffig{fig:prog4_tri-modulus_1d_data_set}.}
  \label{fig:prog4_segmented_regress_result}
\end{figure}

\begin{figure}[tbp]
  \centering
  \includegraphics[scale=0.50]{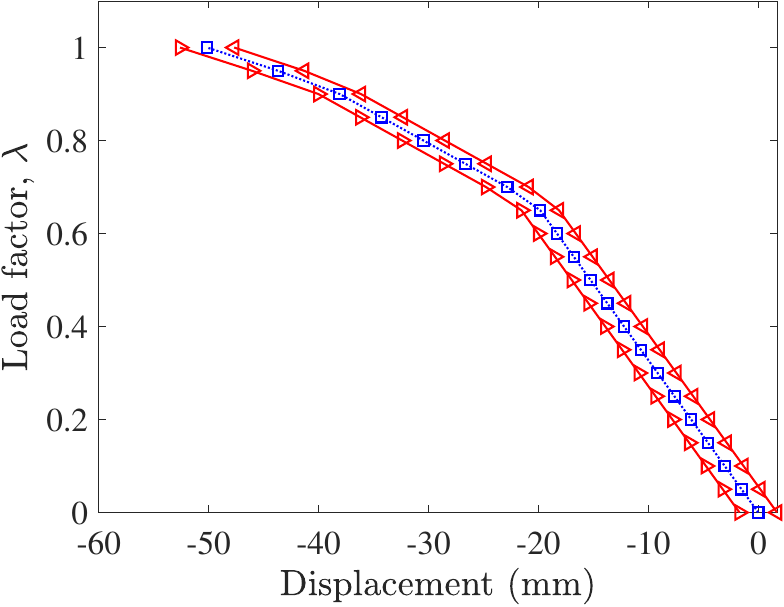}
  \caption{Bounds for the nodal displacement of the frame. 
  ``$\triangleleft$'' and ``$\triangleright$'' denote the upper and 
  lower bounds obtained by the proposed method, respectively; 
  ``$\square$'' denotes the reference solution. }
  \label{fig:prog4_displ_analysis_loop}
\end{figure}

\begin{figure}[tbp]
  \centering
  \includegraphics[scale=0.50]{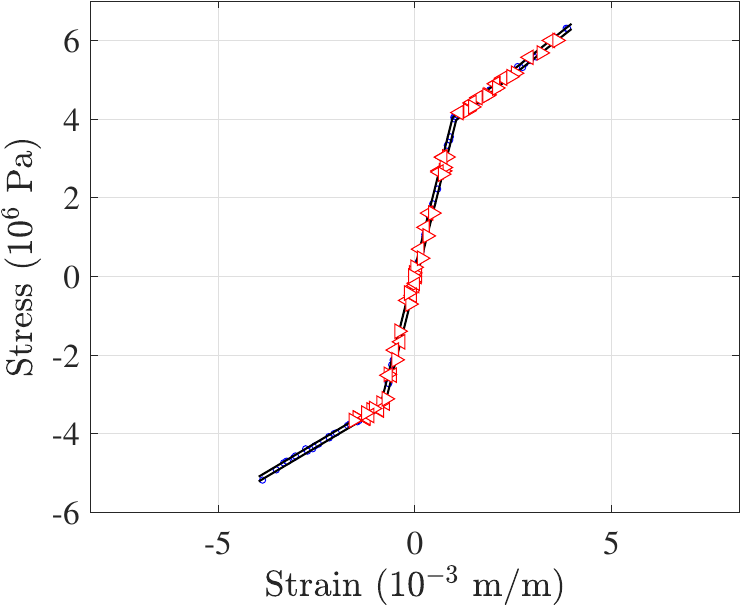}
  \caption{Member stresses--strains corresponding to the obtained 
  solutions at $\lambda=1$ in \reffig{fig:prog4_displ_analysis_loop}. 
  Boundary of the uncertainty set $C(\tau)$ is also depicted. 
  ``$\triangleleft$'' and ``$\triangleright$'' denote the member 
  stresses corresponding to the upper and lower bound solutions, 
  respectively. }
  \label{fig:prog4_displ_analysis_result}
\end{figure}

Consider the planar frame structure illustrated in 
\reffig{fig:prog4_frame}. 
The equilibrium analysis of the frame structure is based on the 
formulation in \citep{Kan16}. 
The top two nodes are fixed, and the frame has $m=12$ members 
and $d=18$ degrees of freedom of the nodal displacements. 
As for the external load, we apply a vertical downward force of 
$2.4 \lambda\,\mathrm{kN}$ at the rightmost node, where 
$\lambda \in [0,1]$ is the load factor. 
Each member has a hollow circular cross-section, and 
the cross-sectional area of each member is $1000\,\mathrm{mm^{2}}$. 
The moment of inertia is computed by supposing that the ratio of the 
inertial radius to the external radius is $0.9$. 

\reffig{fig:prog4_tri-modulus_1d_data_set} shows the material data set 
consisting of $r=80$ data points. 
\reffig{fig:prog4_segmented_regress_result} reports the result of the 
segmented least squares with 
$k=5$, $\mu=2.0$, and $M=1.0\times 10^{3}$. 
We see that the optimal solution uses three straight lines to fit 
the data. 

We set the target reliability and confidence level 
to $1-\epsilon=0.9$ and $1-\delta=0.9$, respectively. 
We adopt the vertical displacement of the rightmost node as the 
quantity of interest. 
\reffig{fig:prog4_displ_analysis_loop} reports the interval bounds 
obtained by the proposed method, as well as the reference solutions, 
for different values of the load factor $\lambda$. 
\reffig{fig:prog4_displ_analysis_result} plots 
$(\varepsilon_{e},\sigma_{e})$ $(e=1,\dots,m)$ corresponding to the 
upper and lower bound solutions at $\lambda=1$, as well as 
the boundary of $C(\tau)$.

\subsection{Cable--strut structure example}
\label{sec:ex.cable}

\begin{figure}[tbp]
  \centering
  \includegraphics[scale=0.60]{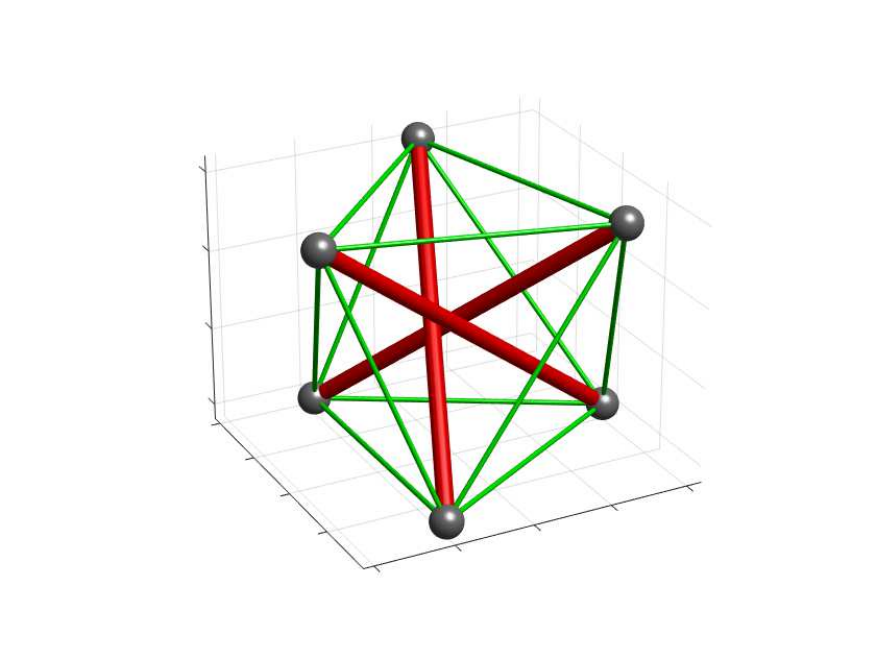}
  \caption{Problem setting of the cable--strut structure example. }
  \label{fig:prog5_cable_strut}
\end{figure}

\begin{figure}[tbp]
  \centering
  \subfloat[]{
  \label{fig:prog5_prog5_tri_no_compress_data_set}
  \includegraphics[scale=0.50]{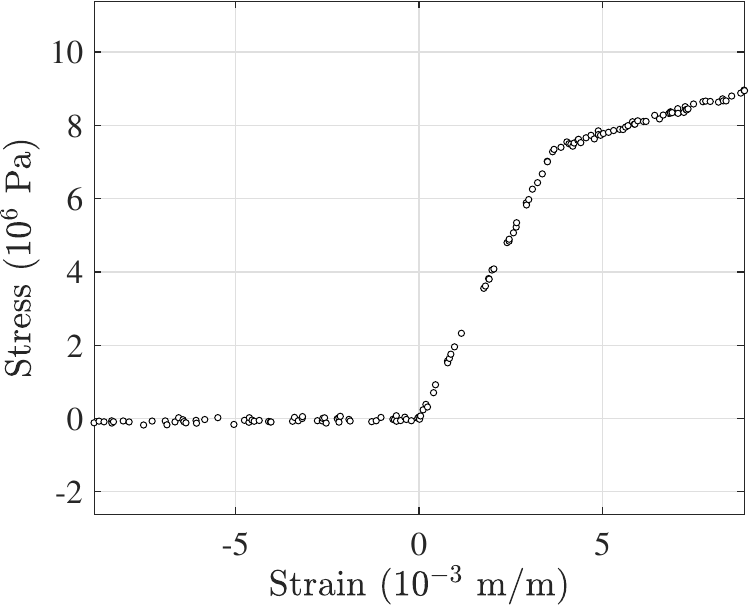}
  }
  \hfill
  \subfloat[]{
  \label{fig:prog5_tri-modulus_1d_data_set}
  \includegraphics[scale=0.50]{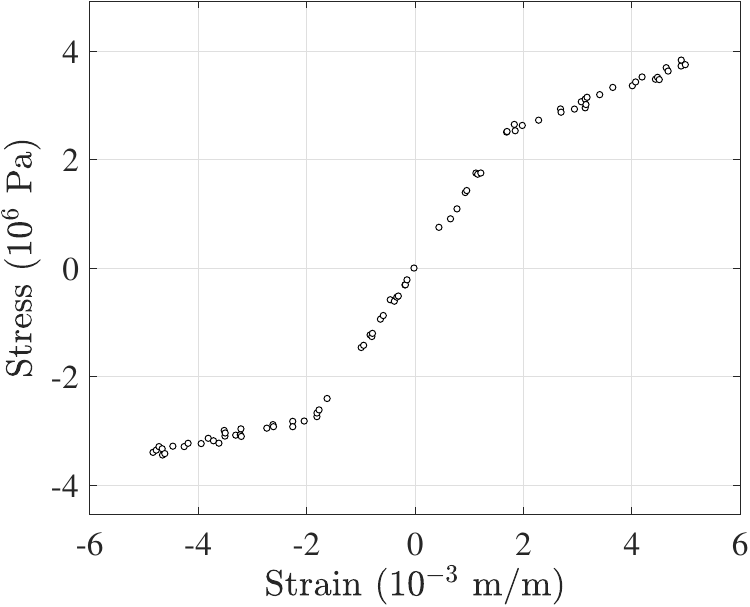}
  }
  \caption[]{Data sets for the cable--strut structure example. 
  \subref{fig:prog5_prog5_tri_no_compress_data_set}~Data set for cables; and 
  \subref{fig:prog5_tri-modulus_1d_data_set}~data set for struts. }
  \label{fig:prog5_data_set}
\end{figure}

\begin{figure}[tbp]
  \centering
  \includegraphics[scale=0.50]{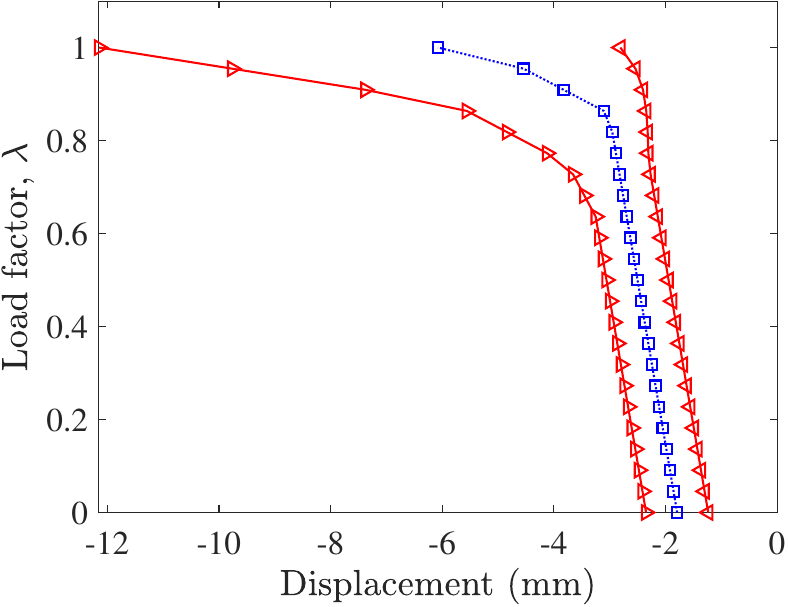}
  \caption{Bounds for the nodal displacement of the cable--strut 
  structure. 
  ``$\triangleleft$'' and ``$\triangleright$'' denote the upper and 
  lower bounds obtained by the proposed method, respectively. 
  ``$\square$'' denote the reference solution. }
  \label{fig:prog5_displ_analysis_loop}
\end{figure}

\begin{figure}[tbp]
  \centering
  \subfloat[]{
  \label{fig:prog5_displ_analysis_cable}
  \includegraphics[scale=0.50]{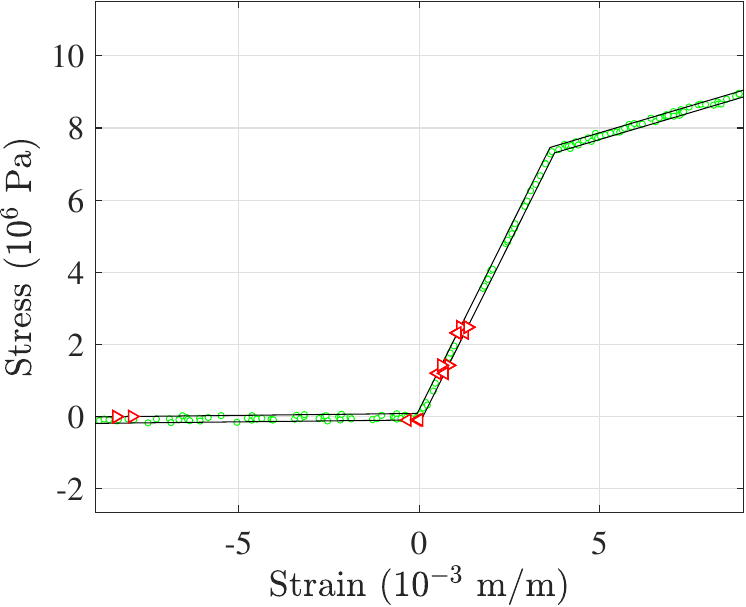}
  }
  \hfill
  \subfloat[]{
  \label{fig:prog5_displ_analysis_strut}
  \includegraphics[scale=0.50]{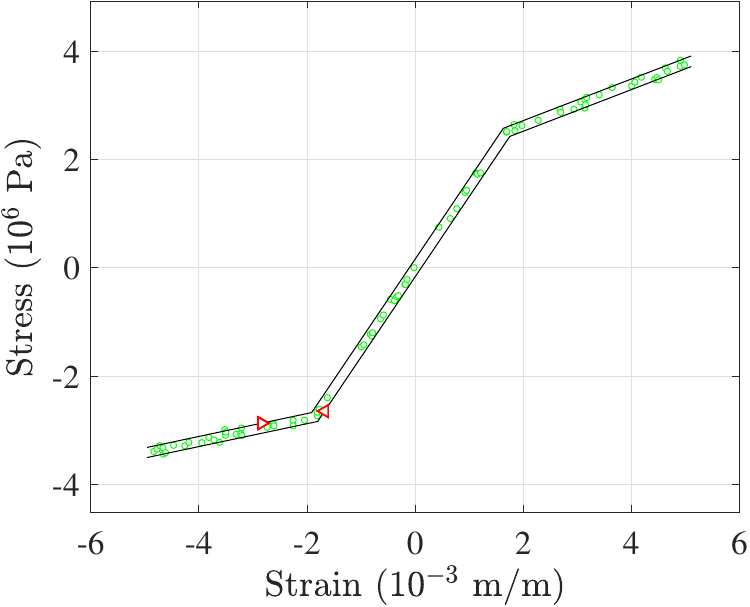}
  }
  \caption[]{Member stresses--strains corresponding to the obtained 
  solutions at $\lambda=1$ in \reffig{fig:prog5_displ_analysis_loop}. 
  Boundary of the uncertainty set $C(\tau)$ is also depicted. 
  ``$\triangleleft$'' and ``$\triangleright$'' denote the member 
  stresses corresponding to the upper and lower bound solutions, 
  respectively. 
  \subref{fig:prog5_displ_analysis_cable}~Cables; and 
  \subref{fig:prog5_displ_analysis_strut}~struts. }
  \label{fig:prog5_displ_analysis}
\end{figure}

\begin{figure}[tbp]
  \centering
  \subfloat[]{
  \label{fig:prog5_displ_analysis_loop_time}
  \includegraphics[scale=0.50]{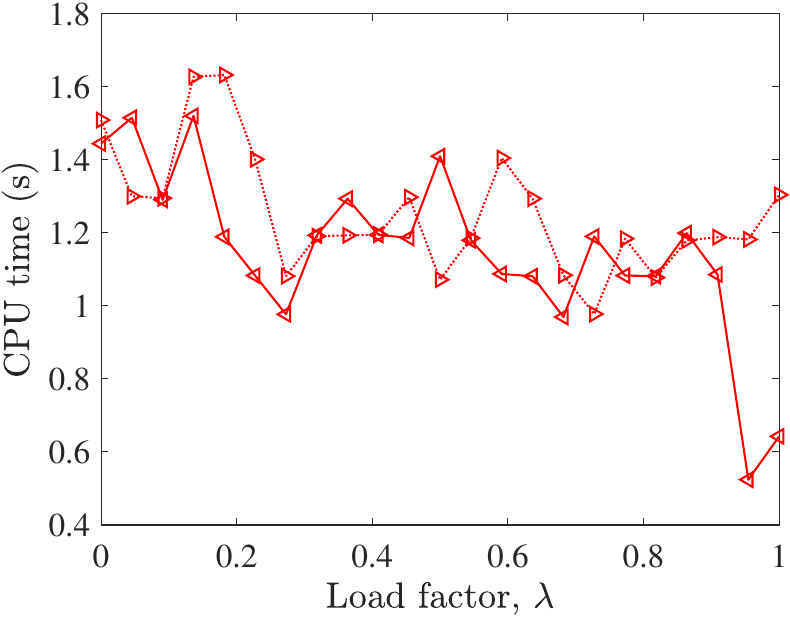}
  }
  \hfill
  \subfloat[]{
  \label{fig:prog5_displ_analysis_loop_node}
  \includegraphics[scale=0.50]{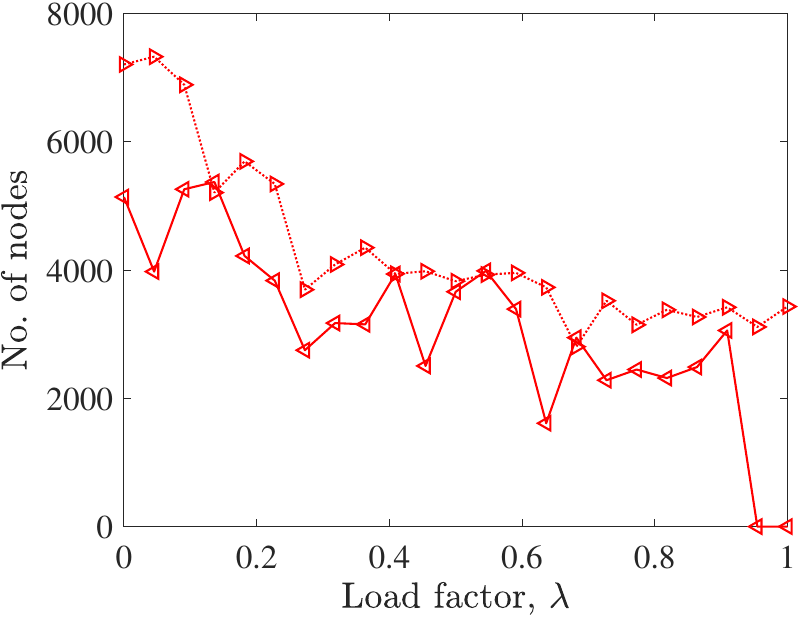}
  }
  \caption[]{
  Computation costs of the MILP problems for the cable--strut structure. 
  ``$\triangleleft$'' and ``$\triangleright$'' correspond to the costs 
  for obtaining the upper and lower bounds in 
  \reffig{fig:prog5_displ_analysis_loop}, respectively. 
  \subref{fig:prog5_displ_analysis_loop_time}~Computation time; and 
  \subref{fig:prog5_displ_analysis_loop_node}~number 
  of enumeration nodes explored by the MILP solver. }
  \label{fig:prog5_displ_analysis_loop_computation}
\end{figure}

Consider the three-dimensional pin-jointed structure 
shown in \reffig{fig:prog5_cable_strut}. 
This structure consists of 12 cables (depicted as the thin lines) and 3 
struts (depicted as the thick lines). 
The bottom and top layers are equilateral triangles consisting of three 
cables, where the length of each cable is $\sqrt{3}\,\mathrm{m}$. 
The bottom layer is in parallel with the top layer with 
$\pi/4\,\mathrm{rad}$ rotated position. 
The distance between these two layers is $1.5\,\mathrm{m}$. 
The length of each strut is $2.4863\,\mathrm{m}$. 
The cross-sectional areas of the cables and struts are 
$500\,\mathrm{mm^{2}}$ and $1000\,\mathrm{mm^{2}}$, respectively. 
To prevent the rigid-body motion, we fix 6 degrees of freedom of the 
displacements of the bottom nodes. 
Accordingly, the number of degrees of freedom of the nodal displacements 
is $d=12$. 
As for the external load, we apply the vertical downward forces 
of $1.1\lambda\,{\mathrm{kN}}$ at the three top nodes, 
where $\lambda \in \Re$ is the load factor. 
As the initial strains, at the nodal location above each cable and each 
strut have strains of $2\times 10^{-3}$ and $-0.4\times 10^{-3}$, 
respectively. 
It is worth noting that the nodal location above is not the 
self-equilibrium shape. 

\reffig{fig:prog5_data_set} collects the material data sets. 
For cables we use the one in 
\reffig{fig:prog5_prog5_tri_no_compress_data_set}, 
which consists of $150$ data points. 
Moreover, for struts we use the one in 
\reffig{fig:prog5_tri-modulus_1d_data_set}, 
which consists of $80$ data points. 
To each data set we apply the segmented least squares with 
$k=5$, $\mu=2.0$, and $M=1.0 \times 10^{3}$. 
In both cases, the optimal solution uses three straight lines (i.e., it 
has two breakpoints). 
By using these results, we compute the reference solutions of the 
equilibrium analysis. 

We set $\epsilon=0.1$ and $\delta=0.1$. 
The quantity of interest is the vertical displacement of a top node. 
\reffig{fig:prog5_displ_analysis_loop} reports the interval bounds 
obtained by the proposed method, as well as the reference solutions. 
As the load factor $\lambda$ increases, the stiffnesses of some cable 
members can be very close to $0$, and hence the interval becomes wider. 
It is confirmed that the reference solution always belongs to the 
obtained interval. 
\reffig{fig:prog5_displ_analysis} plots $(\varepsilon_{e},\sigma_{e})$ 
$(e=1,\dots,m)$ corresponding to the upper and lower bound solutions 
at $\lambda=1$, as well as the boundary of $C(\tau)$. 
Particularly, we can observe in \reffig{fig:prog5_displ_analysis_cable} 
that some cables undergo large compression deformations at the lower 
bound solution. 

\reffig{fig:prog5_displ_analysis_loop_computation} reports the 
computation costs of solving the MILP problems for obtaining the 
solutions shown in \reffig{fig:prog5_displ_analysis_loop}. 
\reffig{fig:prog5_displ_analysis_loop_time} and 
\ref{fig:prog5_displ_analysis_loop_node} shows the computation time and 
the number of enumeration nodes of CPLEX, respectively. 
We can observe that the computation cost is irrelevant to the 
stiffnesses of the cable members.

\section{Conclusions}
\label{sec:conclusions}

This paper has presented an optimization-based method, within the 
framework of data-driven computational elasticity, for computing a 
bound of the structural response considering the uncertainty in the 
material data set. 
The method ensures that, at least the specified confidence level, the 
probability that the structural response is within the obtained bound is 
no smaller than the target reliability. 
This guarantee is provided by a fundamental property of the order 
statistics. 
This means that the proposed method is free from modeling of the 
probabilistic distribution of the constitutive law, and hence 
this method can be viewed as a purely data-driven approach. 

The method developed in this paper is a natural extension 
of the previous work~\citep{Kan23}, which has a drawback that 
its solution overestimates the structural response drastically when the 
stress--strain relation in the given material data set is not 
approximately linear. 
To deal with such a data, this paper has presented a method consisting 
of two steps: We first apply the segmented least squares to the data set; 
then with constructing the uncertainty set based on the order 
statistics we solve optimization problems to find the structural bound. 
In this paper we have shown that the optimization problems in both steps 
can be formulated as mixed-integer convex optimization problems, 
which can be solved globally. 
The guarantee of the global optimality in the latter step is 
particularly crucial, because a bound corresponding to 
a local optimal solution underestimates the structural response in 
general. 
Besides this method, the global optimality is guaranteed in 
the approaches in \citep{Kan23} and \citep{GDLT21}, which aim at handling 
a almost linear stress--strain relation. 
In contrast, the method in \citep{HLDTG23} proposed for dealing with 
nonlinearity lacks guarantee of the global optimality. 

Numerical examples have been demonstrated for three skeletal structures: 
A truss, frame, and cable--strut structure. 
We have compared the interval bounds obtained by the proposed 
method with the reference solutions to confirm validity of the method. 
Since the optimization problems in these examples were solved within 
a few minutes by a standard solver, the computation cost required by 
the proposed method is practically acceptable. 
In contrast, if the given data set fits well a nonlinear 
smooth curve, rather than a piecewise-affine function, then the number 
of 0-1 variables in the proposed formulations increases, which can 
possibly yield drastic increase of computation cost. 
It is worth noting that, besides the structural types in the numerical 
examples presented in this paper, 
we can apply the proposed method to space frames, provided 
that a data set of pairs of shear stresses and shear strains is 
also given. 
An extension to general continua remains as future work.

\paragraph{Acknowledgments}

This work is partially supported by 
JST CREST Grant Number JPMJCR1911, Japan, 
and 
JSPS KAKENHI JP21K04351 and 24K07747. 

\paragraph{Conflict of interest}

The author declares that there is no conflict of interest.


\begin{thebibliography}{99}
\bibitem[\protect\citeauthoryear{Carrara {\em et al.\/}}{2020}]{CDSO20}
  {P.~Carrara, L.~De Lorenzis, L.~Stainier, M.~Ortiz}:
  {Data-driven fracture mechanics}.
  {\em Computer Methods in Applied Mechanics and Engineering},
  \textbf{372}, 113390 (2020).

\bibitem[\protect\citeauthoryear{Ciftci and Hackl}{2022}]{CH22}
  {K.~Ciftci, K.~Hackl}:
  {Model-free data-driven simulation of inelastic materials using 
    structured data sets, tangent space information and transition rules}.
  {\em Computational Mechanics},
  \textbf{70}, 425--435 (2022).

\bibitem[\protect\citeauthoryear{Cl\'{e}ment {\em et al.\/}}{2012}]{CSY12}
  {A.~Cl\'{e}ment, C.~Soize, J.~Yvonnet}:
  {Computational nonlinear stochastic homogenization 
    using a nonconcurrent multiscale approach 
    for hyperelastic heterogeneous microstructures analysis}.
  {\em International Journal for Numerical Methods in Engineering},
  \textbf{91}, 799--824 (2012).

\bibitem[\protect\citeauthoryear{Dal\'{e}mat {\em et al.\/}}{2019}]{DCLV19}
  {M.~Dal\'{e}mat, M.~Coret, A.~Leygue, E.~Verron}:
  {Measuring stress field without constitutive equation}.
  {\em Mechanics of Materials},
  \textbf{136}, 103087 (2019).

\bibitem[\protect\citeauthoryear{Dandin {\em et al.\/}}{2024}]{DLS24}
  {H.~Dandin, A.~Leygue, L.~Stainier}:
  {Graph-based representation of history-dependent material response
    in the data-driven computational mechanics framework}.
  {\em Computer Methods in Applied Mechanics and Engineering},
  \textbf{419}, 116694 (2024).


\bibitem[\protect\citeauthoryear{Eggersmann {\em et al.\/}}{2019}]{EKRSO19}
  {R.~Eggersmann, T.~Kirchdoerfer, S.~Reese, L.~Stainier, M.~Ortiz}:
  {Model-free data-driven inelasticity}.
  {\em Computer Methods in Applied Mechanics and Engineering},
  \textbf{350}, 81--99 (2019).

\bibitem[\protect\citeauthoryear{Ghanem {\em et al.\/}}{2017}]{GHO17}
  {R.~Ghanem, D.~Higdon, H.~Owhadi (eds.)}:
  {\em Handbook of Uncertainty Quantification}.
  Springer International Publishing, Cham (2017).

\bibitem[\protect\citeauthoryear{Guo {\em et al.\/}}{2021}]{GDLT21}
  {X.~Guo, Z.~Du, C.~Liu, S.~Tang}:
  {A new uncertainty analysis-based framework for data-driven 
    computational mechanics}.
  {\em Journal of Applied Mechanics},
  \textbf{88}, 111003 (2021).

\bibitem[\protect\citeauthoryear{Guo {\em et al.\/}}{2023}]{GDWMZSYTG23}
  {Y.~Guo, Z.~Du, L.~Wang, W.~Meng, T.~Zhang, R.~Su, D.~Yang, S.~Tang, X.~Guo}:
  {Data-driven topology optimization (DDTO) for three-dimensional 
    continuum structures}.
  {\em Structural and Multidisciplinary Optimization},
  \textbf{66}, 104 (2023).

\bibitem[\protect\citeauthoryear{Hao {\em et al.\/}}{2022}]{HYYZWW22}
  {P.~Hao, H.~Yang, H.~Yang, Y.~Zhang, Y.~Wang, B.~Wang}:
  {A sequential single-loop reliability optimization and confidence 
    analysis method}.
  {\em Computer Methods in Applied Mechanics and Engineering},
  \textbf{399}, 115400 (2022).

\bibitem[\protect\citeauthoryear{He and Chen}{2020}]{HC20}
  {Q.~He, J.-S.~Chen}:
  {A physics-constrained data-driven approach based on 
    locally convex reconstruction for noisy database}.
  {\em Computer Methods in Applied Mechanics and Engineering},
  \textbf{363}, 112791 (2020).

\bibitem[\protect\citeauthoryear{He {\em et al.\/}}{2021}]{HLLC21}
  {Q.~He, D.~W.~Laurence, C.-H.~Lee, J.-S.~Chen}:
  {Manifold learning based data-driven modeling for soft biological tissues}.
  {\em Journal of Biomechanics},
  \textbf{117}, 110124 (2021).

\bibitem[\protect\citeauthoryear{Huang {\em et al.\/}}{2023}]{HLDTG23}
  {M.~Huang, C.~Liu, Z.~Du, S.~Tang, X.~Guo}:
  {A sequential linear programming (SLP) approach for uncertainty 
    analysis-based data-driven computational mechanics}.
  {\em Computational Mechanics},
  to appear. 
  DOI: \url{10.1007/s00466-023-02395-8}.

\bibitem[\protect\citeauthoryear{Iba\~{n}ez {\em et al.\/}}{2018}]{IAcAGCC18}
  {R.~Iba\~{n}ez, E.~Abisset-Chavanne, J.~V.~Aguado, D.~Gonzalez, E.~Cueto, F.~Chinesta}:
  {A manifold learning approach to data-driven computational 
    elasticity and inelasticity}.
  {\em Archives of Computational Methods in Engineering},
  \textbf{25}, 47--57 (2018).

\bibitem[\protect\citeauthoryear{Iba\~{n}ez {\em et al.\/}}{2017}]{IBAAcCLC17}
  {R.~Iba\~{n}ez, D.~Borzacchiello, J.~V.~Aguado, 
    E.~Abisset-Chavanne, E.~Cueto, P.~Ladeveze, F.~Chinesta}:
  {Data-driven non-linear elasticity: constitutive manifold 
    construction and problem discretization}.
  {\em Computational Mechanics},
  \textbf{60}, 813--826 (2017).

\bibitem[\protect\citeauthoryear{IBM ILOG}{2023}]{cplex}
  {IBM ILOG}:
  {\em IBM ILOG CPLEX Optimization Studio}.
  \url{https://www.ibm.com/products/ilog-cplex-optimization-studio/}
  (Accessed November 2023).

\bibitem[\protect\citeauthoryear{Ito {\em et al.\/}}{2018}]{IKK18}
  {M.~Ito, N.~H.~Kim, N.~Kogiso}:
  {Conservative reliability index for epistemic uncertainty 
    in reliability-based design optimization}.
  {\em Structural and Multidisciplinary Optimization},
  \textbf{57}, 1919--1935 (2018).

\bibitem[\protect\citeauthoryear{Jung {\em et al.\/}}{2020}]{JCDL20}
  {Y.~Jung, H.~Cho, Z.~Duan, I.~Lee}:
  {Determination of sample size for input variables in RBDO 
    through bi-objective confidence-based design optimization 
    under input model uncertainty}.
  {\em Structural and Multidisciplinary Optimization},
  \textbf{61}, 253--266 (2020).

\bibitem[\protect\citeauthoryear{Jung {\em et al.\/}}{2019}]{JCL19}
  {Y.~Jung, H.~Cho, I.~Lee}:
  {Reliability measure approach for confidence-based design optimization 
    under insufficient input data}.
  {\em Structural and Multidisciplinary Optimization},
  \textbf{60}, 1967--1982 (2019).

\bibitem[\protect\citeauthoryear{Kanno}{2016}]{Kan16}
  {Y.~Kanno}:
  {Mixed-integer second-order cone programming for global 
    optimization of compliance of frame structure 
    with discrete design variables}.
  {\em Structural and Multidisciplinary Optimization},
  \textbf{54}, 301--316 (2016).

\bibitem[\protect\citeauthoryear{Kanno}{2019}]{Kan19}
  {Y.~Kanno}:
  {A data-driven approach to non-parametric reliability-based 
    design optimization of structures with uncertain load}.
  {\em Structural and Multidisciplinary Optimization},
  \textbf{60}, 83--97 (2019).

\bibitem[\protect\citeauthoryear{Kanno}{2019}]{Kan19MIP}
  {Y.~Kanno}:
  {Mixed-integer programming formulation of a data-driven solver 
      in computational elasticity}.
  {\em Optimization Letters},
  \textbf{13}, 1505--1514 (2019).

\bibitem[\protect\citeauthoryear{Kanno}{2020}]{Kan20}
  {Y.~Kanno}:
  {Dimensionality reduction enhances data-driven 
    reliability-based design optimizer}.
  {\em Journal of Advanced Mechanical Design, Systems, and Manufacturing},
  \textbf{14}, 19-00200 (2020).

\bibitem[\protect\citeauthoryear{Kanno}{2020}]{Kan20three}
  {Y.~Kanno}:
  {On three concepts in robust design optimization: 
    absolute robustness, relative robustness, and less variance}.
  {\em Structural and Multidisciplinary Optimization},
  \textbf{62}, 979--1000 (2020).

\bibitem[\protect\citeauthoryear{Kanno}{2021}]{Kan21}
  {Y.~Kanno}:
  {A kernel method for learning constitutive relation 
    in data-driven computational elasticity}.
  {\em Japan Journal of Industrial and Applied Mathematics},
  \textbf{38}, 39--77 (2021).

\bibitem[\protect\citeauthoryear{Kanno}{2021}]{Kan21TAML}
  {Y.~Kanno}:
  {Alternating minimization for data-driven computational elasticity from 
    experimental data: kernel method for learning constitutive manifold}.
  {\em Theoretical and Applied Mechanics Letters},
  \textbf{11}, 100289 (2021).

\bibitem[\protect\citeauthoryear{Kanno}{2023}]{Kan23}
  {Y.~Kanno}:
  {Computation-with-confidence for static elasticity: 
    data-driven approach with order statistics}.
  {\em Zeitschrift f\"{u}r Angewandte Mathematik und Mechanik},
  \textbf{103}, e202100482 (2023).

\bibitem[\protect\citeauthoryear{Karapiperis {\em et al.\/}}{2021}]{KSOA21}
  {K.~Karapiperis, L.~Stainier, M.~Ortiz, J.~E.~Andrade}:
  {Data-driven multiscale modeling in mechanics}.
  {\em Journal of the Mechanics and Physics of Solids},
  \textbf{147}, 104239 (2021).

\bibitem[\protect\citeauthoryear{Kirchdoerfer and Ortiz}{2016}]{KO16}
  {T.~Kirchdoerfer, M.~Ortiz}:
  {Data-driven computational mechanics}.
  {\em Computer Methods in Applied Mechanics and Engineering},
  \textbf{304}, 81--101 (2016).

\bibitem[\protect\citeauthoryear{Kleinberg and Tardos}{2006}]{KT06}
  {J.~Kleinberg, \'{E}.~Tardos}:
  {\em Algorithm Design}.
  Pearson Education, Boston (2006).

\bibitem[\protect\citeauthoryear{Leygue {\em et al.\/}}{2018}]{LCRSV18}
  {A.~Leygue, M.~Coret, J.~R\'{e}thor\'{e}, L.~Stainier, E.~Verron}:
  {Data-based derivation of material response}.
  {\em Computer Methods in Applied Mechanics and Engineering},
  \textbf{331}, 184--196 (2018).


\bibitem[\protect\citeauthoryear{Luo and Paal}{2023}]{LP23}
  {H.~Luo, S.~G.~Paal}:
  {A novel outlier-insensitive local support vector machine 
    for robust data-driven forecasting in engineering}.
  {\em Engineering with Computers},
  \textbf{39}, 3671--3689 (2023).

\bibitem[\protect\citeauthoryear{Moon {\em et al.\/}}{2018}]{MCCGLG18}
  {M.-Y.~Moon, H.~Cho, K.~K.~Choi, N.~Gaul, D.~Lamb, D.~Gorsich}:
  {Confidence-based reliability assessment considering 
    limited numbers of both input and output test data}.
  {\em Structural and Multidisciplinary Optimization},
  \textbf{57}, 2027--2043 (2018).

\bibitem[\protect\citeauthoryear{Mora-Mac\'{i}as {\em et al.\/}}{2020}]{MARDDDS20}
  {J.~Mora-Mac\'{i}as, J.~Ayensa-Jim\'{e}nez, E.~Reina-Romo, 
    M.~H.~Doweidar, J.~Dom\'{i}nguez, M.~Doblar\'{e}, J.~A.~Sanz-Herrera}:
  {A multiscale data-driven approach for bone tissue biomechanics}.
  {\em Computer Methods in Applied Mechanics and Engineering},
  \textbf{368}, 113136, (2020).

\bibitem[\protect\citeauthoryear{Nguyen and Keip}{2018}]{NK18}
  {L.~T.~K.~Nguyen, M.-A.~Keip}:
  {A data-driven approach to nonlinear elasticity}.
  {\em Computers and Structures},
  \textbf{194}, 97--115 (2018).

\bibitem[\protect\citeauthoryear{Nguyen {\em et al.\/}}{2020}]{NRK20}
  {L.~T.~K.~Nguyen, M.~Rambausek, M.-A.~Keip}:
  {Variational framework for distance-minimizing method 
    in data-driven computational mechanics}.
  {\em Computer Methods in Applied Mechanics and Engineering},
  \textbf{365}, 112898 (2020).

\bibitem[\protect\citeauthoryear{Pham {\em et al.\/}}{2023}]{PBG23}
  {D.~K.~N.~Pham, N.~Blal, A.~Gravouil}:
  {Tangent space data driven framework for elasto-plastic material behaviors}.
  {\em Finite Elements in Analysis and Design},
  \textbf{216}, 103895 (2023).


\bibitem[\protect\citeauthoryear{Poelstra {\em et al.\/}}{2023}]{PBS23}
  {K.~Poelstra, T.~Bartel, B.~Schweizer}:
  {A data-driven framework for evolutionary problems in solid mechanics}.
  {\em Zeitschrift f\"{u}r Angewandte Mathematik und Mechanik},
  \textbf{103}, e202100538 (2023).

\bibitem[\protect\citeauthoryear{Prume {\em et al.\/}}{2023}]{PRO23}
  {E.~Prume, S.~Reese, M.~Ortiz}:
  {Model-free data-driven inference in computational mechanics}.
  {\em Computer Methods in Applied Mechanics and Engineering},
  \textbf{403}, 115704 (2023).



\bibitem[\protect\citeauthoryear{Stainier {\em et al.\/}}{2019}]{SLO19}
  {L.~Stainier, A.~Leygue, M.~Ortiz}:
  {Model-free data-driven methods in mechanics: material data 
    identification and solvers}.
  {\em Computational Mechanics},
  \textbf{64}, 381--393 (2019).

\bibitem[\protect\citeauthoryear{Su {\em et al.\/}}{2023}]{SJC23}
  {T.-H.~Su, J.~G.~Jean, C.-S.~Chen}:
  {Model-free data-driven identification algorithm enhanced 
    by local manifold learning}.
  {\em Computational Mechanics volume},
  \textbf{71}, 637--655 (2023).

\bibitem[\protect\citeauthoryear{Tang {\em et al.\/}}{2020}]{TLQYSMLG20}
  {S.~Tang, Y.~Li, H.~Qiu, H.~Yang, S.~Saha, S.~Mojumder, W.~K.~Liu, X.~Guo}:
  {MAP123-EP: A mechanistic-based data-driven approach 
    for numerical elastoplastic analysis}.
  {\em Computer Methods in Applied Mechanics and Engineering},
  \textbf{364}, 112955 (2020).

\bibitem[\protect\citeauthoryear{Tang {\em et al.\/}}{2021}]{TYQFLG21}
  {S.~Tang, H.~Yang, H.~Qiu, M.~Fleming, W.~K~Liu, X.~Guo}:
  {MAP123-EPF: A mechanistic-based data-driven approach 
    for numerical elastoplastic modeling at finite strain}.
  {\em Computer Methods in Applied Mechanics and Engineering},
  \textbf{373}, 113484 (2021).

\bibitem[\protect\citeauthoryear{Tang {\em et al.\/}}{2019}]{TZYLLG19}
  {S.~Tang, G.~Zhang, H.~Yang, Y.~Li, W.~K.~Liu, X.~Guo}:
  {MAP123: a data-driven approach to use 1D data for 3D nonlinear 
    elastic materials modeling}.
  {\em Computer Methods in Applied Mechanics and Engineering},
  \textbf{357}, 112587 (2019).

\bibitem[\protect\citeauthoryear{Temizer and Wriggers}{2007}]{TW07}
  {\.{I}.~Temizer, P.~Wriggers}:
  {An adaptive method for homogenization in orthotropic nonlinear elasticity}.
  {\em Computer Methods in Applied Mechanics and Engineering},
  \textbf{196}, 3409--3423 (2007).

\bibitem[\protect\citeauthoryear{Terada {\em et al.\/}}{2013}]{TKHIY13}
  {K.~Terada, J.~Kato, N.~Hirayama, T.~Inugai, K.~Yamamoto}:
  {A method of two-scale analysis with micro-macro decoupling scheme: 
    application to hyperelastic composite materials}.
  {\em Computational Mechanics},
  \textbf{52}, 1199--1219 (2013).

\bibitem[\protect\citeauthoryear{Wang {\em et al.\/}}{2020}]{WHYWG20}
  {Y.~Wang, P.~Hao, H.~Yang, B.~Wang, Q.~Gao}:
  {A confidence-based reliability optimization with single loop 
    strategy and second-order reliability method}.
  {\em Computer Methods in Applied Mechanics and Engineering},
  \textbf{372}, 113436 (2020).

\bibitem[\protect\citeauthoryear{Watanabe and Terada}{2010}]{WT10}
  {I.~Watanabe, K.~Terada}:
  {A method of predicting macroscopic yield strength of polycrystalline 
    metals subjected to plastic forming by micro-macro de-coupling scheme}.
  {\em International Journal of Mechanical Sciences},
  \textbf{52}, 343--355 (2010).

\bibitem[\protect\citeauthoryear{Xu {\em et al.\/}}{2020}]{XYYHGBZBH20}
  {R.~Xu, J.~Yang, W.~Yan, Q.~Huang, G.~Giunta, S.~Belouettar, 
    H.~Zahrouni, T.~Ben Zineb, H.~Hu}:
  {Data-driven multiscale finite element method: 
    from concurrence to separation},
  {\em Computer Methods in Applied Mechanics and Engineering},
  \textbf{363}, 112893 (2020).

\end{thebibliography}
\end{document}